# MAXIMUM LIKELIHOOD ESTIMATION FOR $\alpha$-STABLE AUTOREGRESSIVE PROCESSES

By Beth Andrews[1], Matthew Calder[2] and Richard A. Davis[1,2,3]

*Northwestern University, PHZ Capital Partners and Columbia University*

We consider maximum likelihood estimation for both causal and noncausal autoregressive time series processes with non-Gaussian $\alpha$-stable noise. A nondegenerate limiting distribution is given for maximum likelihood estimators of the parameters of the autoregressive model equation and the parameters of the stable noise distribution. The estimators for the autoregressive parameters are $n^{1/\alpha}$-consistent and converge in distribution to the maximizer of a random function. The form of this limiting distribution is intractable, but the shape of the distribution for these estimators can be examined using the bootstrap procedure. The bootstrap is asymptotically valid under general conditions. The estimators for the parameters of the stable noise distribution have the traditional $n^{1/2}$ rate of convergence and are asymptotically normal. The behavior of the estimators for finite samples is studied via simulation, and we use maximum likelihood estimation to fit a noncausal autoregressive model to the natural logarithms of volumes of Wal-Mart stock traded daily on the New York Stock Exchange.

**1. Introduction.** Many observed time series processes appear "spiky" due to the occasional appearance of observations particularly large in absolute value. Non-Gaussian $\alpha$-stable distributions, which have regularly varying or "heavy" tail probabilities $(\mathrm{P}(|X| > x) \sim (constant)x^{-\alpha}, x > 0, 0 < \alpha < 2)$, are often used to model these series. Processes exhibiting non-Gaussian stable behavior have appeared, for example, in economics and finance (Embrechts, Klüppelberg and Mikosch [18], McCulloch [25] and Mittnik and

Received November 2007; revised July 2008.

[1]Supported in part by NSF Grant DMS-03-08109.

[2]Supported in part by NSF Grant DMS-95-04596.

[3]Supported in part by NSF Grant DMS-07-43459.

Rachev [28]), signal processing (Nikias and Shao [29]) and teletraffic engineering (Resnick [32]).

The focus of this paper is maximum likelihood (ML) estimation for the parameters of autoregressive (AR) time series processes with non-Gaussian stable noise. Specific applications for heavy-tailed AR models include fitting network interarrival times (Resnick [32]), sea surface temperatures (Gallagher [20]) and stock market log-returns (Ling [24]). Causality (all roots of the AR polynomial are outside the unit circle in the complex plane) is a common assumption in the time series literature since causal and non-causal models are indistinguishable in the case of Gaussian noise. However, noncausal AR models are identifiable in the case of non-Gaussian noise, and these models are frequently used in deconvolution problems (Blass and Halsey [3], Chien, Yang and Chi [10], Donoho [16] and Scargle [36]) and have also appeared for modeling stock market trading volume data (Breidt, Davis and Trindade [5]). We, therefore, consider parameter estimation for both causal and noncausal AR models. We assume the parameters of the AR model equation and the parameters of the stable noise distribution are unknown, and we maximize the likelihood function with respect to all parameters. Since most stable density functions do not have a closed-form expression, the likelihood function is evaluated by inversion of the stable characteristic function. We show that ML estimators of the AR parameters are $n^{1/\alpha}$-consistent ($n$ represents sample size) and converge in distribution to the maximizer of a random function. The form of this limiting distribution is intractable, but the shape of the distribution for these estimators can be examined using the bootstrap procedure. We show the bootstrap procedure is asymptotically valid provided the bootstrap sample size $m_n \to \infty$ with $m_n/n \to 0$ as $n \to \infty$. ML estimators of the parameters of the stable noise distribution are $n^{1/2}$-consistent, asymptotically independent of the AR estimators and have a multivariate normal limiting distribution.

Parameter estimation for causal, heavy-tailed AR processes has already been considered in the literature (Davis and Resnick [14], least squares estimators; Davis [11] and Davis, Knight and Liu [12], least absolute deviations and other $M$-estimators; Mikosch, Gadrich, Klüppelberg and Adler [27], Whittle estimators; Ling [24], weighted least absolute deviations estimators). The weighted least absolute deviations estimators for causal AR parameters are $n^{1/2}$-consistent, and the least squares and Whittle estimators are $(n/\ln n)^{1/\alpha}$-consistent, while the unweighted least absolute deviations estimators have the same faster rate of convergence as ML estimators, $n^{1/\alpha}$. Least absolute deviations and ML estimators have different limiting distributions, however, and simulation results in Calder and Davis [8] show that ML estimates (obtained using the stable likelihood) tend to be more efficient than least absolute deviations estimates, even when the AR process has regularly varying tail probabilities but is not stable. Theory has not yet been developed for the distribution of AR parameter estimators when the process is noncausal and heavy-tailed.

In Section 2, we discuss properties of AR processes with non-Gaussian stable noise and give an approximate log-likelihood for the model parameters. In Section 3, we give a nondegenerate limiting distribution for ML estimators, show that the bootstrap procedure can be used to approximate the distribution for AR parameter estimators, and discuss confidence interval calculation for the model parameters. Proofs of the Lemmas used to establish the results of Section 3 can be found in the Appendix. We study the behavior of the estimators for finite samples via simulation in Section 4.1 and, in Section 4.2, use ML estimation to fit a noncausal AR model to the natural logarithms of volumes of Wal-Mart stock traded daily on the New York Stock Exchange. A causal AR model is inadequate for these log-volumes since causal AR residuals appear dependent. The noncausal residuals appear i.i.d. (independent and identically distributed) stable, and so the fitted noncausal AR model appears much more suitable for the series.

## 2. Preliminaries.
Let $\{X_t\}$ be the AR process which satisfies the difference equations

$$\phi_0(B)X_t = Z_t, \tag{2.1}$$

where the AR polynomial $\phi_0(z) := 1 - \phi_{01}z - \cdots - \phi_{0p}z^p \neq 0$ for $|z| = 1$, $B$ is the backshift operator ($B^k X_t = X_{t-k}$, $k = 0, \pm 1, \pm 2, \ldots$), and $\{Z_t\}$ is an i.i.d. sequence of random variables. Because $\phi_0(z) \neq 0$ for $|z| = 1$, the Laurent series expansion of $1/\phi_0(z)$, $1/\phi_0(z) = \sum_{j=-\infty}^{\infty} \psi_j z^j$, exists on some annulus $\{z : a^{-1} < |z| < a\}$, $a > 1$, and the unique, strictly stationary solution to (2.1) is given by $X_t = \sum_{j=-\infty}^{\infty} \psi_j Z_{t-j}$ (see Brockwell and Davis [6], Chapter 3). Note that if $\phi_0(z) \neq 0$ for $|z| \leq 1$, then $\psi_j = 0$ for $j < 0$ and so $\{X_t\}$ is said to be *causal* since $X_t = \sum_{j=0}^{\infty} \psi_j Z_{t-j}$, a function of only the past and present $\{Z_t\}$. On the other hand, if $\phi_0(z) \neq 0$ for $|z| \geq 1$, then $X_t = \sum_{j=0}^{\infty} \psi_{-j} Z_{t+j}$, and $\{X_t\}$ is said to be a *purely noncausal* process. In the purely noncausal case, the coefficients $\{\psi_j\}$ satisfy $(1 - \phi_{01}z - \cdots - \phi_{0p}z^p)(\psi_0 + \psi_{-1}z^{-1} + \cdots) = 1$, which, if $\phi_{0p} \neq 0$, implies that $\psi_0 = \psi_{-1} = \cdots = \psi_{1-p} = 0$ and $\psi_{-p} = -\phi_{0p}^{-1}$. To express $\phi_0(z)$ as the product of causal and purely noncausal polynomials, suppose

$$\phi_0(z) = (1 - \theta_{01}z - \cdots - \theta_{0r_0}z^{r_0})(1 - \theta_{0,r_0+1}z - \cdots - \theta_{0,r_0+s_0}z^{s_0}), \tag{2.2}$$

where $r_0 + s_0 = p$, $\theta_0^\dagger(z) := 1 - \theta_{01}z - \cdots - \theta_{0r_0}z^{r_0} \neq 0$ for $|z| \leq 1$, and $\theta_0^*(z) := 1 - \theta_{0,r_0+1}z - \cdots - \theta_{0,r_0+s_0}z^{s_0} \neq 0$ for $|z| \geq 1$. Hence, $\theta_0^\dagger(z)$ is a causal polynomial and $\theta_0^*(z)$ is a purely noncausal polynomial. So that $\phi_0(z)$ has a unique representation as the product of causal and purely noncausal polynomials $\theta_0^\dagger(z)$ and $\theta_0^*(z)$, if the true order of the polynomial $\phi_0(z)$ is less than $p$ (if $\phi_{0p} = 0$), we further suppose that $\theta_{0,r_0+s_0} \neq 0$ when $s_0 > 0$. Therefore, if the true order of the AR polynomial $\phi_0(z)$ is less than $p = r_0 + s_0$, then the true order of $\theta_0^\dagger(z)$ is less than $r_0$, but the order of $\theta_0^*(z)$ is $s_0$.

We assume throughout that the i.i.d. noise $\{Z_t\}$ have a univariate stable distribution with exponent $\alpha_0 \in (0,2)$, parameter of symmetry $|\beta_0| < 1$, scale parameter $0 < \sigma_0 < \infty$, and location parameter $\mu_0 \in \mathbb{R}$. Let $\boldsymbol{\tau}_0 = (\alpha_0, \beta_0, \sigma_0, \mu_0)'$. By definition, nondegenerate, i.i.d. random variables $\{S_t\}$ have a stable distribution if there exist positive constants $\{a_n\}$ and constants $\{b_n\}$ such that $a_n(S_1 + \cdots + S_n) + b_n \overset{\mathcal{L}}{=} S_1$ for all $n$. In general, stable distributions are indexed by an exponent $\alpha \in (0,2]$, a parameter of symmetry $|\beta| \leq 1$, a scale parameter $0 < \sigma < \infty$ and a location parameter $\mu \in \mathbb{R}$. Hence, $\boldsymbol{\tau}_0$ is in the interior of the stable parameter space. If $\beta = 0$, the stable distribution is symmetric about $\mu$, and, if $\alpha = 1$ and $\beta = 0$, the symmetric distribution is Cauchy. When $\alpha = 2$, the stable distribution is Gaussian with mean $\mu$ and standard deviation $\sqrt{2}\sigma$. Other properties of stable distributions can be found in Feller [19], Gnedenko and Kolmogorov [21], Samorodnitsky and Taqqu [35] and Zolotarev [38].

Since the stable noise distribution has exponent $\alpha_0 < 2$,

$$\lim_{x \to \infty} x^{\alpha_0} \mathrm{P}(|Z_t| > x) = \tilde{c}(\alpha_0)\sigma_0^{\alpha_0},$$

(2.3)

$$\text{with } \tilde{c}(\alpha) := \left( \int_0^\infty t^{-\alpha} \sin(t)\, dt \right)^{-1}$$

(Samorodnitsky and Taqqu [35], Property 1.2.15). Following Properties 1.2.1 and 1.2.3 in Samorodnitsky and Taqqu [35], $X_t = \sum_{j=-\infty}^{\infty} \psi_j Z_{t-j}$ also has a stable distribution with exponent $\alpha_0$ and, hence, the tail probabilities for the AR process $\{X_t\}$ are also regularly varying with exponent $\alpha_0$. It follows that $\mathrm{E}|X_t|^\delta < \infty$ for all $\delta \in [0, \alpha_0)$ and $\mathrm{E}|X_t|^\delta = \infty$ for all $\delta \geq \alpha_0$.

The characteristic function for $Z_t$ is

(2.4)
$$\varphi_0(s) := \mathrm{E}\{\exp(isZ_t)\}$$
$$= \begin{cases} \exp\left\{ -\sigma_0^{\alpha_0}|s|^{\alpha_0}\left[1 + i\beta_0(\mathrm{sign}\,s)\tan\left(\frac{\pi\alpha_0}{2}\right)(|\sigma_0 s|^{1-\alpha_0} - 1)\right] \\ \qquad\qquad\qquad\qquad\qquad\qquad\qquad\qquad\qquad + i\mu_0 s \right\}, \\ \qquad \text{if } \alpha_0 \neq 1, \\ \exp\left\{ -\sigma_0|s|\left[1 + i\beta_0\frac{2}{\pi}(\mathrm{sign}\,s)\ln(\sigma_0|s|)\right] + i\mu_0 s \right\}, \\ \qquad \text{if } \alpha_0 = 1, \end{cases}$$

and so the density function for the noise can be expressed as $f(z; \boldsymbol{\tau}_0) = (2\pi)^{-1} \int_{-\infty}^{\infty} \exp(-izs)\varphi_0(s)\, ds$. No general, closed-form expression is known for $f$, however; although, computational formulas exist that can be used to evaluate $f$ (see, e.g., McCulloch [26] and Nolan [30]). It can be shown that $f(z; \boldsymbol{\tau}_0) = \sigma_0^{-1} f(\sigma_0^{-1}(z - \mu_0); (\alpha_0, \beta_0, 1, 0)')$, $f(\cdot; (\alpha_0, \beta_0, 1, 0)')$ is unimodal on $\mathbb{R}$ (Yamazato [37]), and $f(z; (\alpha, \beta, 1, 0)')$ is infinitely differentiable with respect to $(z, \alpha, \beta)$ on $\mathbb{R} \times (0, 2) \times (-1, 1)$. There are alternative parameterizations for the stable characteristic function $\varphi_0$ (see, e.g., Zolotarev [38]), but we are using (2.4) so that the noise density function is differentiable with respect to not only $z$ on $\mathbb{R}$ but also $(\alpha, \beta, \sigma, \mu)'$ on $(0, 2) \times (-1, 1) \times (0, \infty) \times (-\infty, \infty)$. From asymptotic expansions in Du-Mouchel [17], if $\boldsymbol{\Omega}_\delta := \{\boldsymbol{\tau} = (\alpha, \beta, \sigma, \mu)' : \|\boldsymbol{\tau} - \boldsymbol{\tau}_0\| < \delta\}$, then for $\delta > 0$ sufficiently small we have the following bounds for the partial and mixed partial derivatives of $\ln f(z; \boldsymbol{\tau})$ as $|z| \to \infty$:

$$\text{(2.5)} \quad \bullet \quad \sup_{\boldsymbol{\Omega}_\delta} \left| \frac{\partial^2 \ln f(z; \boldsymbol{\tau})}{\partial z^2} \right| + \sup_{\boldsymbol{\Omega}_\delta} \left| \frac{\partial^2 \ln f(z; \boldsymbol{\tau})}{\partial z \, \partial \mu} \right| + \sup_{\boldsymbol{\Omega}_\delta} \left| \frac{\partial^2 \ln f(z; \boldsymbol{\tau})}{\partial \mu^2} \right|$$
$$= O(|z|^{-2}),$$

$$\text{(2.6)} \quad \bullet \quad \sup_{\boldsymbol{\Omega}_\delta} \left| \frac{\partial \ln f(z; \boldsymbol{\tau})}{\partial z} \right| + \sup_{\boldsymbol{\Omega}_\delta} \left| \frac{\partial \ln f(z; \boldsymbol{\tau})}{\partial \mu} \right| + \sup_{\boldsymbol{\Omega}_\delta} \left| \frac{\partial^2 \ln f(z; \boldsymbol{\tau})}{\partial z \, \partial \beta} \right|$$
$$+ \sup_{\boldsymbol{\Omega}_\delta} \left| \frac{\partial^2 \ln f(z; \boldsymbol{\tau})}{\partial z \, \partial \sigma} \right| + \sup_{\boldsymbol{\Omega}_\delta} \left| \frac{\partial^2 \ln f(z; \boldsymbol{\tau})}{\partial \beta \, \partial \mu} \right|$$
$$+ \sup_{\boldsymbol{\Omega}_\delta} \left| \frac{\partial^2 \ln f(z; \boldsymbol{\tau})}{\partial \sigma \, \partial \mu} \right| = O(|z|^{-1}),$$

$$\text{(2.7)} \quad \bullet \quad \sup_{\boldsymbol{\Omega}_\delta} \left| \frac{\partial^2 \ln f(z; \boldsymbol{\tau})}{\partial z \, \partial \alpha} \right| + \sup_{\boldsymbol{\Omega}_\delta} \left| \frac{\partial^2 \ln f(z; \boldsymbol{\tau})}{\partial \alpha \, \partial \mu} \right| = O(|z|^{-1} \ln |z|),$$

$$\text{(2.8)} \quad \bullet \quad \sup_{\boldsymbol{\Omega}_\delta} \left| \frac{\partial \ln f(z; \boldsymbol{\tau})}{\partial \beta} \right| + \sup_{\boldsymbol{\Omega}_\delta} \left| \frac{\partial \ln f(z; \boldsymbol{\tau})}{\partial \sigma} \right| + \sup_{\boldsymbol{\Omega}_\delta} \left| \frac{\partial^2 \ln f(z; \boldsymbol{\tau})}{\partial \beta^2} \right|$$
$$+ \sup_{\boldsymbol{\Omega}_\delta} \left| \frac{\partial^2 \ln f(z; \boldsymbol{\tau})}{\partial \beta \, \partial \sigma} \right| + \sup_{\boldsymbol{\Omega}_\delta} \left| \frac{\partial^2 \ln f(z; \boldsymbol{\tau})}{\partial \sigma^2} \right| = O(1),$$

$$\text{(2.9)} \quad \bullet \quad \sup_{\boldsymbol{\Omega}_\delta} \left| \frac{\partial \ln f(z; \boldsymbol{\tau})}{\partial \alpha} \right| + \sup_{\boldsymbol{\Omega}_\delta} \left| \frac{\partial^2 \ln f(z; \boldsymbol{\tau})}{\partial \alpha \, \partial \beta} \right| + \sup_{\boldsymbol{\Omega}_\delta} \left| \frac{\partial^2 \ln f(z; \boldsymbol{\tau})}{\partial \alpha \, \partial \sigma} \right|$$
$$= O(\ln |z|),$$

$$\text{(2.10)} \quad \bullet \quad \sup_{\boldsymbol{\Omega}_\delta} \left| \frac{\partial^2 \ln f(z; \boldsymbol{\tau})}{\partial \alpha^2} \right| = O([\ln |z|]^2).$$

From (2.1) and (2.2), $Z_t = (1 - \theta_{01} B - \cdots - \theta_{0r_0} B^{r_0})(1 - \theta_{0,r_0+1} B - \cdots - \theta_{0,r_0+s_0} B^{s_0}) X_t$. Therefore, for arbitrary autoregressive polynomials $\theta^\dagger(z) = 1 - \theta_1 z - \cdots - \theta_r z^r$ and $\theta^*(z) = 1 - \theta_{r+1} z - \cdots - \theta_{r+s} z^s$, with $r + s = p$, $\theta^{\dagger}(z) \neq 0$ for $|z| \leq 1$, $\theta^{*}(z) \neq 0$ for $|z| \geq 1$, and $\theta_{r+s} \neq 0$ when $s > 0$, we define

$$(2.11) \quad Z_t(\boldsymbol{\theta}, s) = (1 - \theta_1 B - \cdots - \theta_r B^r)(1 - \theta_{r+1} B - \cdots - \theta_{r+s} B^s) X_t,$$

where $\boldsymbol{\theta} := (\theta_1, \ldots, \theta_p)'$. Let $\boldsymbol{\theta}_0 = (\theta_{01}, \ldots, \theta_{0p})'$ denote the true parameter vector and note that $\{Z_t(\boldsymbol{\theta}_0, s_0)\} = \{Z_t\}$. Now, let $\boldsymbol{\eta} = (\eta_1, \ldots, \eta_{p+4})' = (\theta_1, \ldots, \theta_p, \alpha, \beta, \sigma, \mu)' = (\boldsymbol{\theta}', \boldsymbol{\tau}')'$, and let $\boldsymbol{\eta}_0 = (\eta_{01}, \ldots, \eta_{0,p+4})' = (\boldsymbol{\theta}_0', \boldsymbol{\tau}_0')'$. From Breidt et al. [4], given a realization $\{X_t\}_{t=1}^n$ from (2.1), the log-likelihood of $\boldsymbol{\eta}$ can be approximated by the conditional log-likelihood

$$(2.12) \quad \mathcal{L}(\boldsymbol{\eta}, s) = \sum_{t=p+1}^{n} \left[ \ln f(Z_t(\boldsymbol{\theta}, s); \boldsymbol{\tau}) + \ln |\theta_p| I\{s > 0\} \right],$$

where $\{Z_t(\boldsymbol{\theta}, s)\}_{t=p+1}^n$ is computed using (2.11) and $I\{\cdot\}$ represents the indicator function (see [4] for the derivation of $\mathcal{L}$). Given $\{X_t\}_{t=1}^n$ and fixed $p$, we can estimate $s_0$, the order of noncausality for the AR model (2.1), and $\boldsymbol{\eta}_0$ by maximizing $\mathcal{L}$ with respect to both $s$ and $\boldsymbol{\eta}$. If the function $g$ is defined so that

$$g(\boldsymbol{\theta}, s) = [g_j(\boldsymbol{\theta}, s)]_{j=1}^{p},$$

$$(2.13) \quad g_j(\boldsymbol{\theta}, s) = \begin{cases} \theta_j - \displaystyle\sum_{k=1}^{j} \theta_{j-k} \theta_{p-s+k}, & j = 1, \ldots, p-s, \\ -\displaystyle\sum_{k=j-p+s}^{j} \theta_{j-k} \theta_{p-s+k}, & j = p-s+1, \ldots, p, \end{cases}$$

with $\theta_0 = -1$ and $\theta_k = 0$ whenever $k \notin \{0, \ldots, p\}$, then an estimate of $\boldsymbol{\phi}_0 := (\phi_{01}, \ldots, \phi_{0p})'$ can be obtained using the MLEs of $s_0$ and $\boldsymbol{\theta}_0$ and the fact that $\boldsymbol{\phi}_0 = g(\boldsymbol{\theta}_0, s_0)$. A similar ML approach is considered in [4] for lighter-tailed AR processes.

**3. Asymptotic results.** In this section, we obtain limiting results for maximizers of the log-likelihood $\mathcal{L}$. But first, we need to introduce some notation and define a random function $W(\cdot)$. The ML estimators of $\boldsymbol{\theta}_0$ converge in distribution to the maximizer of $W(\cdot)$.

Suppose the Laurent series expansions for $1/\theta_0^{\dagger}(z) = 1/(1 - \theta_{01} z - \cdots - \theta_{0r_0} z^{r_0})$ and $1/\theta_0^{*}(z) = 1/(1 - \theta_{0,r_0+1} z - \cdots - \theta_{0,r_0+s_0} z^{s_0})$ are given by $1/\theta_0^{\dagger}(z) = \sum_{j=0}^{\infty} \pi_j z^j$ and $1/\theta_0^{*}(z) = \sum_{j=s_0}^{\infty} \chi_j z^{-j}$. From (2.11),

$$(3.1) \quad \frac{\partial Z_t(\boldsymbol{\theta}, s)}{\partial \theta_j} = \begin{cases} -\theta^{*}(B) X_{t-j}, & j = 1, \ldots, r, \\ -\theta^{\dagger}(B) X_{t+r-j}, & j = r+1, \ldots, p, \end{cases}$$

and so, for $\mathbf{u} = (u_1, \ldots, u_p)' \in \mathbb{R}^p$,

$$
\begin{aligned}
\mathbf{u}'\frac{\partial Z_t(\boldsymbol{\theta}_0, s_0)}{\partial \boldsymbol{\theta}} &= -u_1\theta_0^*(B)X_{t-1} - \cdots - u_{r_0}\theta_0^*(B)X_{t-r_0} - u_{r_0+1}\theta_0^\dagger(B)X_{t-1} \\
&\quad - \cdots - u_p\theta_0^\dagger(B)X_{t-s_0} \\
&= -u_1(1/\theta_0^\dagger(B))Z_{t-1} - \cdots - u_{r_0}(1/\theta_0^\dagger(B))Z_{t-r_0} \\
&\quad - u_{r_0+1}(1/\theta_0^*(B))Z_{t-1} - \cdots - u_p(1/\theta_0^*(B))Z_{t-s_0} \\
&= -u_1\sum_{j=0}^{\infty}\pi_j Z_{t-1-j} - \cdots - u_{r_0}\sum_{j=0}^{\infty}\pi_j Z_{t-r_0-j} \\
&\quad - u_{r_0+1}\sum_{j=s_0}^{\infty}\chi_j Z_{t-1+j} - \cdots - u_p\sum_{j=s_0}^{\infty}\chi_j Z_{t-s_0+j}.
\end{aligned}
$$

Therefore, if

$$
\tag{3.2} \sum_{j=-\infty}^{\infty} c_j(\mathbf{u})Z_{t-j} := \mathbf{u}'\frac{\partial Z_t(\boldsymbol{\theta}_0, s_0)}{\partial \boldsymbol{\theta}},
$$

then $c_0(\mathbf{u}) = -u_p\chi_{s_0}I\{s_0 > 0\} = u_p\theta_{0p}^{-1}I\{s_0 > 0\}$, $c_1(\mathbf{u}) = -u_1\pi_0 I\{r_0 > 0\} = -u_1 I\{r_0 > 0\}$, $c_{-1}(\mathbf{u}) = -u_p\chi_2 I\{s_0 = 1\} - (u_{p-1}\chi_{s_0} + u_p\chi_{s_0+1})I\{s_0 > 1\}$, and so on. Since $\{\pi_j\}_{j=0}^{\infty}$ and $\{\chi_j\}_{j=s_0}^{\infty}$ decay at geometric rates (Brockwell and Davis [6], Chapter 3), for any $\mathbf{u} \in \mathbb{R}^p$, there exist constants $C(\mathbf{u}) > 0$ and $0 < D(\mathbf{u}) < 1$ such that

$$
\tag{3.3} |c_j(\mathbf{u})| \leq C(\mathbf{u})[D(\mathbf{u})]^{|j|} \qquad \forall j \in \{\ldots, -1, 0, 1, \ldots\}.
$$

We now define the function

$$
\tag{3.4}
\begin{aligned}
W(\mathbf{u}) = \sum_{k=1}^{\infty}\sum_{j\neq 0}\{&\ln f(Z_{k,j} + [\tilde{c}(\alpha_0)]^{1/\alpha_0}\sigma_0 c_j(\mathbf{u})\delta_k\Gamma_k^{-1/\alpha_0}; \boldsymbol{\tau}_0) \\
&- \ln f(Z_{k,j}; \boldsymbol{\tau}_0)\},
\end{aligned}
$$

where

- $\{Z_{k,j}\}_{k,j}$ is an i.i.d. sequence with $Z_{k,j} \stackrel{\mathcal{L}}{=} Z_1$,
- $\tilde{c}(\cdot)$ was defined in (2.3),
- $\{\delta_k\}$ is i.i.d. with $\mathrm{P}(\delta_k = 1) = (1 + \beta_0)/2$ and $\mathrm{P}(\delta_k = -1) = 1 - (1 + \beta_0)/2$,
- $\Gamma_k = E_1 + \cdots + E_k$, where $\{E_k\}$ is an i.i.d. series of exponential random variables with mean one, and
- $\{Z_{k,j}\}$, $\{\delta_k\}$ and $\{E_k\}$ are mutually independent.

Note that $(1 + \beta_0)/2 = \lim_{x\to\infty}[\mathrm{P}(Z_1 > x)/\mathrm{P}(|Z_1| > x)]$ (Samorodnitsky and Taqqu [35], Property 1.2.15). Some properties of $W(\cdot)$ are given in the following theorem.

THEOREM 3.1. *With probability one, the function $W(\mathbf{u})$ defined in (3.4) is finite for all $\mathbf{u} \in \mathbb{R}^p$ and has a unique maximum.*

PROOF. Let $\mathbf{u} \in \mathbb{R}^p$ and observe that

$$W(\mathbf{u}) = \sum_{k=1}^{\infty} \sum_{j \neq 0} [\tilde{c}(\alpha_0)]^{1/\alpha_0} \sigma_0 c_j(\mathbf{u}) \delta_k \Gamma_k^{-1/\alpha_0} \frac{\partial \ln f(Z_{k,j}(\mathbf{u}); \boldsymbol{\tau}_0)}{\partial z}$$

$$= \sum_{k=1}^{\infty} \sum_{j \neq 0} [\tilde{c}(\alpha_0)]^{1/\alpha_0} \sigma_0 c_j(\mathbf{u}) \delta_k$$

$$\times \Gamma_k^{-1/\alpha_0} \left[ \frac{\partial \ln f(Z_{k,j}(\mathbf{u}); \boldsymbol{\tau}_0)}{\partial z} - \frac{\partial \ln f(Z_{k,j}; \boldsymbol{\tau}_0)}{\partial z} \right]$$

$$+ \sum_{k=1}^{\infty} \sum_{j \neq 0} [\tilde{c}(\alpha_0)]^{1/\alpha_0} \sigma_0 c_j(\mathbf{u}) \delta_k (\Gamma_k^{-1/\alpha_0} - k^{-1/\alpha_0}) \frac{\partial \ln f(Z_{k,j}; \boldsymbol{\tau}_0)}{\partial z}$$

$$+ \sum_{k=1}^{\infty} \sum_{j \neq 0} [\tilde{c}(\alpha_0)]^{1/\alpha_0} \sigma_0 c_j(\mathbf{u}) \delta_k k^{-1/\alpha_0} \frac{\partial \ln f(Z_{k,j}; \boldsymbol{\tau}_0)}{\partial z},$$

where $\tilde{Z}_{k,j}(\mathbf{u})$ lies between $Z_{k,j}$ and $Z_{k,j} + [\tilde{c}(\alpha_0)]^{1/\alpha_0} \sigma_0 c_j(\mathbf{u}) \delta_k \Gamma_k^{-1/\alpha_0}$. Since $0 < [\tilde{c}(\alpha_0)]^{1/\alpha_0} \sigma_0 < \infty$, by Lemmas A.1–A.3 in the Appendix, $|W(\mathbf{u})| < \infty$ almost surely. It can be shown similarly that $\sup_{\|\mathbf{u}\| \leq T} |W(\mathbf{u})| < \infty$ almost surely for any $T \in (0, \infty)$ and, therefore, $\mathrm{P}(\bigcap_{T=1}^{\infty} \{\sup_{\|\mathbf{u}\| \leq T} |W(\mathbf{u})| < \infty\}) = 1$.

Since $f(\cdot; \boldsymbol{\tau}_0)$ is unimodal and differentiable on $\mathbb{R}$, with positive probability, $\ln f(Z_1 + \cdot; \boldsymbol{\tau}_0)$ is strictly concave in a neighborhood of zero, and so, by Remark 2 in Davis, Knight and Liu [12], $W(\cdot)$ has a unique maximum almost surely. □

We now give nondegenerate limiting distributions for ML estimators of $\boldsymbol{\eta}_0 = (\boldsymbol{\theta}_0', \boldsymbol{\tau}_0')' = (\theta_{01}, \ldots, \theta_{0p}, \alpha_0, \beta_0, \sigma_0, \mu_0)'$ and estimators of the AR parameters $\boldsymbol{\phi}_0 = (\phi_{01}, \ldots, \phi_{0p})'$ in (2.1).

THEOREM 3.2. *There exists a sequence of maximizers $\hat{\boldsymbol{\eta}}_{\mathrm{ML}} = (\hat{\boldsymbol{\theta}}_{\mathrm{ML}}', \hat{\boldsymbol{\tau}}_{\mathrm{ML}}')'$ of $\mathcal{L}(\cdot, s_0)$ in (2.12) such that, as $n \to \infty$,*

$$(3.5) \quad n^{1/\alpha_0}(\hat{\boldsymbol{\theta}}_{\mathrm{ML}} - \boldsymbol{\theta}_0) \xrightarrow{\mathcal{L}} \boldsymbol{\xi} \quad \text{and} \quad n^{1/2}(\hat{\boldsymbol{\tau}}_{\mathrm{ML}} - \boldsymbol{\tau}_0) \xrightarrow{\mathcal{L}} \mathbf{Y} \sim N(\mathbf{0}, \mathbf{I}^{-1}(\boldsymbol{\tau}_0)),$$

*where $\boldsymbol{\xi}$ is the unique maximizer of $W(\cdot)$, $\boldsymbol{\xi}$ and $\mathbf{Y}$ are independent, and $\mathbf{I}(\boldsymbol{\tau}) := -[E\{\partial^2 \ln f(Z_1; \boldsymbol{\tau})/(\partial \tau_i \partial \tau_j)\}]_{i,j \in \{1, \ldots, 4\}}$. In addition, if $\hat{\boldsymbol{\phi}}_{\mathrm{ML}} := g(\hat{\boldsymbol{\theta}}_{\mathrm{ML}}, s_0)$, with $g$ as defined in (2.13), then*

$$(3.6) \quad n^{1/\alpha_0}(\hat{\boldsymbol{\phi}}_{\mathrm{ML}} - \boldsymbol{\phi}_0) \xrightarrow{\mathcal{L}} \boldsymbol{\Sigma}(\boldsymbol{\theta}_0)\boldsymbol{\xi},$$

*where*

(3.7)
$$\boldsymbol{\Sigma}(\boldsymbol{\theta}) := \begin{bmatrix} \dfrac{\partial g_1(\boldsymbol{\theta}, s_0)}{\partial \theta_1} & \cdots & \dfrac{\partial g_1(\boldsymbol{\theta}, s_0)}{\partial \theta_p} \\ \vdots & \ddots & \vdots \\ \dfrac{\partial g_p(\boldsymbol{\theta}, s_0)}{\partial \theta_1} & \cdots & \dfrac{\partial g_p(\boldsymbol{\theta}, s_0)}{\partial \theta_p} \end{bmatrix}$$

*and $g_1, \ldots, g_p$ were also defined in (2.13).*

Since $\boldsymbol{\tau}_0$ is in the interior of the stable parameter space, given i.i.d. observations $\{Z_t\}_{t=1}^n$, ML estimators of $\boldsymbol{\tau}_0$ are asymptotically Gaussian with mean $\boldsymbol{\tau}_0$ and covariance matrix $\mathbf{I}^{-1}(\boldsymbol{\tau}_0)/n$ (see DuMouchel [17]). The estimators $\hat{\boldsymbol{\tau}}_{\mathrm{ML}}$, therefore, have the same limiting distribution as ML estimators in the case of observed i.i.d. noise. Nolan [31] lists values of $\mathbf{I}^{-1}(\cdot)$ for different parameter values.

For $\mathbf{u} \in \mathbb{R}^p$ and $\mathbf{v} \in \mathbb{R}^4$, let $W_n(\mathbf{u}, \mathbf{v}) = \mathcal{L}(\boldsymbol{\eta}_0 + (n^{-1/\alpha_0}\mathbf{u}', n^{-1/2}\mathbf{v}')', s_0) - \mathcal{L}(\boldsymbol{\eta}_0, s_0)$, and note that maximizing $\mathcal{L}(\boldsymbol{\eta}, s_0)$ with respect to $\boldsymbol{\eta}$ is equivalent to maximizing $W_n(\mathbf{u}, \mathbf{v})$ with respect to $\mathbf{u}$ and $\mathbf{v}$ if $\mathbf{u} = n^{1/\alpha_0}(\boldsymbol{\theta} - \boldsymbol{\theta}_0)$ and $\mathbf{v} = n^{1/2}(\boldsymbol{\tau} - \boldsymbol{\tau}_0)$. We give a functional convergence result for $W_n$ in the following theorem, and then use it to prove Theorem 3.2.

THEOREM 3.3. *As $n \to \infty$, $W_n(\mathbf{u}, \mathbf{v}) \xrightarrow{\mathcal{L}} W(\mathbf{u}) + \mathbf{v}'\mathbf{N} - 2^{-1}\mathbf{v}'\mathbf{I}(\boldsymbol{\tau}_0)\mathbf{v}$ on $C(\mathbb{R}^{p+4})$, where $\mathbf{N} \sim N(\mathbf{0}, \mathbf{I}(\boldsymbol{\tau}_0))$ is independent of $W(\cdot)$, and $C(\mathbb{R}^{p+4})$ represents the space of continuous functions on $\mathbb{R}^{p+4}$, where convergence is equivalent to uniform convergence on every compact subset.*

PROOF. For $\mathbf{u} \in \mathbb{R}^p$ and $\mathbf{v} \in \mathbb{R}^4$, let

$$W_n^*(\mathbf{u}, \mathbf{v}) = \sum_{t=p+1}^n \left\{ \ln f\left(Z_t + n^{-1/\alpha_0} \sum_{j \neq 0} c_j(\mathbf{u}) Z_{t-j}; \boldsymbol{\tau}_0\right) - \ln f(Z_t; \boldsymbol{\tau}_0) \right\}$$
$$+ \frac{\mathbf{v}'}{\sqrt{n}} \sum_{t=p+1}^n \frac{\partial \ln f(Z_t; \boldsymbol{\tau}_0)}{\partial \boldsymbol{\tau}}.$$

Since

$$W_n(\mathbf{u}, \mathbf{v}) - W_n^*(\mathbf{u}, \mathbf{v})$$
$$= \sum_{t=p+1}^n \ln f\left(Z_t\left(\boldsymbol{\theta}_0 + \frac{\mathbf{u}}{n^{1/\alpha_0}}, s_0\right); \boldsymbol{\tau}_0 + \frac{\mathbf{v}}{\sqrt{n}}\right)$$
$$- \sum_{t=p+1}^n \ln f\left(Z_t + n^{-1/\alpha_0} \sum_{j \neq 0} c_j(\mathbf{u}) Z_{t-j}; \boldsymbol{\tau}_0\right)$$

$$- \frac{\mathbf{v}'}{\sqrt{n}} \sum_{t=p+1}^{n} \frac{\partial \ln f(Z_t; \boldsymbol{\tau}_0)}{\partial \boldsymbol{\tau}} + (n-p) \ln \left| \frac{\theta_{0p} + n^{-1/\alpha_0} u_p}{\theta_{0p}} \right| I\{s_0 > 0\},$$

$W_n(\mathbf{u}, \mathbf{v}) - W_n^*(\mathbf{u}, \mathbf{v}) + 2^{-1} \mathbf{v}' \mathbf{I}(\boldsymbol{\tau}_0) \mathbf{v} = o_p(1)$ on $C(\mathbb{R}^{p+4})$ by Lemmas A.4–A.7. So, the proof is complete if $W_n^*(\mathbf{u}, \mathbf{v}) \overset{\mathcal{L}}{\to} W(\mathbf{u}) + \mathbf{v}' \mathbf{N}$ on $C(\mathbb{R}^{p+4})$.

For $\mathbf{u} \in \mathbb{R}^p$, let

$$(3.8) \quad W_n^{\dagger}(\mathbf{u}) = \sum_{t=p+1}^{n} \left[ \ln f\left( Z_t + n^{-1/\alpha_0} \sum_{j \neq 0} c_j(\mathbf{u}) Z_{t-j}; \boldsymbol{\tau}_0 \right) - \ln f(Z_t; \boldsymbol{\tau}_0) \right]$$

and, for $\mathbf{v} \in \mathbb{R}^4$, let

$$(3.9) \quad T_n(\mathbf{v}) = \frac{\mathbf{v}'}{\sqrt{n}} \sum_{t=p+1}^{n} \frac{\partial \ln f(Z_t; \boldsymbol{\tau}_0)}{\partial \boldsymbol{\tau}}.$$

By Lemma A.8, for fixed $\mathbf{u}$ and $\mathbf{v}$, $(W_n^{\dagger}(\mathbf{u}), T_n(\mathbf{v}))' \overset{\mathcal{L}}{\to} (W(\mathbf{u}), \mathbf{v}' \mathbf{N})'$ on $\mathbb{R}^2$, with $W(\mathbf{u})$ and $\mathbf{v}' \mathbf{N}$ independent. Consequently, $W_n^*(\mathbf{u}, \mathbf{v}) = W_n^{\dagger}(\mathbf{u}) + T_n(\mathbf{v}) \overset{\mathcal{L}}{\to} W(\mathbf{u}) + \mathbf{v}' \mathbf{N}$ on $\mathbb{R}$. Similarly, it can be shown that the finite dimensional distributions of $W_n^*(\mathbf{u}, \mathbf{v})$ converge to those of $W(\mathbf{u}) + \mathbf{v}' \mathbf{N}$, with $W(\cdot)$ and $\mathbf{N}$ independent. For any compact set $K_1 \subset \mathbb{R}^p$, $\{W_n^{\dagger}(\cdot)\}$ is tight on $C(K_1)$ by Lemma A.12 and, for any compact set $K_2 \subset \mathbb{R}^4$, $\{T_n(\cdot)\}$ is tight on $C(K_2)$ since $T_n(\mathbf{v})$ is linear in $\mathbf{v}$. Therefore, by Theorem 7.1 in Billingsley [2], $W_n^*(\mathbf{u}, \mathbf{v}) = W_n^{\dagger}(\mathbf{u}) + T_n(\mathbf{v}) \overset{\mathcal{L}}{\to} W(\mathbf{u}) + \mathbf{v}' \mathbf{N}$ on $C(\mathbb{R}^{p+4})$. $\square$

PROOF OF THEOREM 3.2. Since $W_n(\mathbf{u}, \mathbf{v}) \overset{\mathcal{L}}{\to} W(\mathbf{u}) + \mathbf{v}' \mathbf{N} - 2^{-1} \mathbf{v}' \mathbf{I}(\boldsymbol{\tau}_0) \mathbf{v}$ on $C(\mathbb{R}^{p+4})$, $\boldsymbol{\xi}$ uniquely maximizes $W(\cdot)$ almost surely, and $\mathbf{Y} = \mathbf{I}^{-1}(\boldsymbol{\tau}_0) \mathbf{N}$ uniquely maximizes $\mathbf{v}' \mathbf{N} - 2^{-1} \mathbf{v}' \mathbf{I}(\boldsymbol{\tau}_0) \mathbf{v}$, from Remark 1 in Davis, Knight and Liu [12], there exists a sequence of maximizers of $W_n(\cdot, \cdot)$ which converges in distribution to $(\boldsymbol{\xi}', \mathbf{Y}')'$. The result (3.5) follows because $\mathcal{L}(\boldsymbol{\eta}, s_0) - \mathcal{L}(\boldsymbol{\eta}_0, s_0) = W_n(n^{1/\alpha_0}(\boldsymbol{\theta} - \boldsymbol{\theta}_0), n^{1/2}(\boldsymbol{\tau} - \boldsymbol{\tau}_0))$. By Theorem 3.3, $\boldsymbol{\xi}$ and $\mathbf{Y}$ are independent.

Using the mean-value theorem,

$$n^{1/\alpha_0}(\hat{\boldsymbol{\phi}}_{\text{ML}} - \boldsymbol{\phi}_0) = n^{1/\alpha_0}(g(\hat{\boldsymbol{\theta}}_{\text{ML}}, s_0) - g(\boldsymbol{\theta}_0, s_0))$$

$$(3.10) \qquad = \begin{bmatrix} \dfrac{\partial g_1(\boldsymbol{\theta}_1^*, s_0)}{\partial \theta_1} & \cdots & \dfrac{\partial g_1(\boldsymbol{\theta}_1^*, s_0)}{\partial \theta_p} \\ \vdots & \ddots & \vdots \\ \dfrac{\partial g_p(\boldsymbol{\theta}_p^*, s_0)}{\partial \theta_1} & \cdots & \dfrac{\partial g_p(\boldsymbol{\theta}_p^*, s_0)}{\partial \theta_p} \end{bmatrix}$$

$$\times n^{1/\alpha_0}(\hat{\boldsymbol{\theta}}_{\text{ML}} - \boldsymbol{\theta}_0),$$

where $\boldsymbol{\theta}_1^*, \ldots, \boldsymbol{\theta}_p^*$ lie between $\hat{\boldsymbol{\theta}}_{\mathrm{ML}}$ and $\boldsymbol{\theta}_0$. Since $\hat{\boldsymbol{\theta}}_{\mathrm{ML}} \xrightarrow{P} \boldsymbol{\theta}_0$ and $\boldsymbol{\Sigma}(\cdot)$ is continuous at $\boldsymbol{\theta}_0$, (3.10) equals $\boldsymbol{\Sigma}(\boldsymbol{\theta}_0) n^{1/\alpha_0}(\hat{\boldsymbol{\theta}}_{\mathrm{ML}} - \boldsymbol{\theta}_0) + o_p(1)$. Therefore, the result (3.6) follows from (3.5). □

Since the forms of the limiting distributions for $\hat{\boldsymbol{\theta}}_{\mathrm{ML}}$ and $\hat{\boldsymbol{\phi}}_{\mathrm{ML}}$ in (3.5) and (3.6) are intractable, we recommend using the bootstrap procedure to examine the distributions for these estimators. Davis and Wu [15] give a bootstrap procedure for examining the distribution of $M$-estimates for the parameters of causal, heavy-tailed AR processes; we consider a similar procedure here. Given observations $\{X_t\}_{t=1}^n$ from (2.1), $\hat{\boldsymbol{\theta}}_{\mathrm{ML}}$ from (3.5), and corresponding residuals $\{Z_t(\hat{\boldsymbol{\theta}}_{\mathrm{ML}}, s_0)\}_{t=p+1}^n$ obtained via (2.11), the procedure is implemented by first generating an i.i.d. sequence $\{Z_t^*\}_{t=1}^{m_n}$ from the empirical distribution for $\{Z_t(\hat{\boldsymbol{\theta}}_{\mathrm{ML}}, s_0)\}_{t=p+1}^n$. A bootstrap replicate $X_1^*, \ldots, X_{m_n}^*$ is then obtained from the fitted AR($p$) model

$$(3.11) \qquad \hat{\theta}_{\mathrm{ML}}^{\dagger}(B) \hat{\theta}_{\mathrm{ML}}^*(B) X_t^* = Z_t^*,$$

where $\hat{\theta}_{\mathrm{ML}}^{\dagger}(z) := 1 - \hat{\theta}_{1,\mathrm{ML}} z - \cdots - \hat{\theta}_{r_0,\mathrm{ML}} z^{r_0}$ and $\hat{\theta}_{\mathrm{ML}}^*(z) := 1 - \hat{\theta}_{r_0+1,\mathrm{ML}} z - \cdots - \hat{\theta}_{r_0+s_0,\mathrm{ML}} z^{s_0}$ (let $Z_t^* = 0$ for $t \notin \{1, \ldots, m_n\}$). Finally, with $Z_t^*(\boldsymbol{\theta}, s) := (1 - \theta_1 B - \cdots - \theta_r B^r)(1 - \theta_{r+1} B - \cdots - \theta_{r+s} B^s) X_t^*$ for $\boldsymbol{\theta} = (\theta_1, \ldots, \theta_p)' \in \mathbb{R}^p$ and $r + s = p$, a bootstrap replicate $\hat{\boldsymbol{\theta}}_{m_n}^*$ of $\hat{\boldsymbol{\theta}}_{\mathrm{ML}}$ can be found by maximizing

$$\mathcal{L}_{m_n}^*(\boldsymbol{\theta}, s_0) := \sum_{t=p+1}^{m_n} \left[ \ln f(Z_t^*(\boldsymbol{\theta}, s_0); \hat{\boldsymbol{\tau}}_{\mathrm{ML}}) + \ln |\theta_p| I\{s_0 > 0\} \right]$$

with respect to $\boldsymbol{\theta}$. The limiting behavior of $\hat{\boldsymbol{\theta}}_{m_n}^*$, along with that of $\hat{\boldsymbol{\phi}}_{m_n}^* := g(\hat{\boldsymbol{\theta}}_{m_n}^*, s_0)$ (a bootstrap replicate of $\hat{\boldsymbol{\phi}}_{\mathrm{ML}}$), is considered in Theorem 3.4. To give a precise statement of the results, we let $\mathcal{M}_p(\mathbb{R}^p)$ represent the space of probability measures on $\mathbb{R}^p$ and we use the metric $d_p$ from Davis and Wu ([15], page 1139) to metrize the topology of weak convergence on $\mathcal{M}_p(\mathbb{R}^p)$. For random elements $Q_n$ and $Q$ of $\mathcal{M}_p(\mathbb{R}^p)$, $Q_n \xrightarrow{P} Q$ if and only if $d_p(Q_n, Q) \xrightarrow{P} 0$ on $\mathbb{R}$, which is equivalent to $\int_{\mathbb{R}^p} h_j \, dQ_n \xrightarrow{P} \int_{\mathbb{R}^p} h_j \, dQ$ on $\mathbb{R}$ for all $j \in \{1, 2, \ldots\}$, where $\{h_j\}_{j=1}^{\infty}$ is a dense sequence of bounded, uniformly continuous functions on $\mathbb{R}^p$. By Theorem 3.4, $\mathrm{P}(m_n^{1/\hat{\alpha}_{\mathrm{ML}}}(\hat{\boldsymbol{\theta}}_{m_n}^* - \hat{\boldsymbol{\theta}}_{\mathrm{ML}}) \in \cdot | X_1, \ldots, X_n)$ converges in probability to $\mathrm{P}(\boldsymbol{\xi} \in \cdot)$ on $\mathcal{M}_p(\mathbb{R}^p)$ [$\boldsymbol{\xi}$ represents the unique maximizer of $W(\cdot)$], and a similar result holds for $m_n^{1/\hat{\alpha}_{\mathrm{ML}}}(\hat{\boldsymbol{\phi}}_{m_n}^* - \hat{\boldsymbol{\phi}}_{\mathrm{ML}})$.

THEOREM 3.4. *If, as $n \to \infty$, $m_n \to \infty$ with $m_n/n \to 0$, then there exists a sequence of maximizers $\hat{\boldsymbol{\theta}}_{m_n}^*$ of $\mathcal{L}_{m_n}^*(\cdot, s_0)$ such that*

$$P(m_n^{1/\hat{\alpha}_{\mathrm{ML}}}(\hat{\boldsymbol{\theta}}_{m_n}^* - \hat{\boldsymbol{\theta}}_{\mathrm{ML}}) \in \cdot | X_1, \ldots, X_n) \xrightarrow{P} P(\boldsymbol{\xi} \in \cdot)$$

on $\mathcal{M}_p(\mathbb{R}^p)$ and, if $\hat{\boldsymbol{\phi}}_{m_n}^* = g(\hat{\boldsymbol{\theta}}_{m_n}^*, s_0)$, then

$$(3.12) \qquad P(m_n^{1/\hat{\alpha}_{\mathrm{ML}}}(\hat{\boldsymbol{\phi}}_{m_n}^* - \tilde{\boldsymbol{\phi}}_{\mathrm{ML}}) \in \cdot | X_1, \ldots, X_n) \overset{P}{\to} P(\boldsymbol{\Sigma}(\boldsymbol{\theta}_0)\boldsymbol{\xi} \in \cdot)$$

on $\mathcal{M}_p(\mathbb{R}^p)$ [$\boldsymbol{\Sigma}(\cdot)$ was defined in (3.7)].

PROOF.  Since $Z_t^*(\boldsymbol{\theta}, s) = (1 - \theta_1 B - \cdots - \theta_r B^r)(1 - \theta_{r+1}B - \cdots - \theta_{r+s}B^s) \times X_t^*$, following (3.1), for $\mathbf{u} = (u_1, \ldots, u_p)' \in \mathbb{R}^p$,

$$
\begin{aligned}
\mathbf{u}'\frac{\partial Z_t^*(\hat{\boldsymbol{\theta}}_{\mathrm{ML}}, s_0)}{\partial \boldsymbol{\theta}} &= -u_1 \hat{\theta}_{\mathrm{ML}}^*(B) X_{t-1}^* - \cdots - u_{r_0} \hat{\theta}_{\mathrm{ML}}^*(B) X_{t-r_0}^* \\
&\quad - u_{r_0+1} \hat{\theta}_{\mathrm{ML}}^\dagger(B) X_{t-1}^* - \cdots - u_p \hat{\theta}_{\mathrm{ML}}^\dagger(B) X_{t-s_0}^* \\
&= -u_1(1/\hat{\theta}_{\mathrm{ML}}^\dagger(B)) Z_{t-1}^* - \cdots - u_{r_0}(1/\hat{\theta}_{\mathrm{ML}}^\dagger(B)) Z_{t-r_0}^* \\
&\quad - u_{r_0+1}(1/\hat{\theta}_{\mathrm{ML}}^*(B)) Z_{t-1}^* - \cdots - u_p(1/\hat{\theta}_{\mathrm{ML}}^*(B)) Z_{t-s_0}^*.
\end{aligned}
$$

We define the sequence $\{\hat{c}_j(\mathbf{u})\}_{j=-\infty}^{\infty}$ so that

$$(3.13) \qquad \sum_{j=-\infty}^{\infty} \hat{c}_j(\mathbf{u}) Z_{t-j}^* = \mathbf{u}'\frac{\partial Z_t^*(\hat{\boldsymbol{\theta}}_{\mathrm{ML}}, s_0)}{\partial \boldsymbol{\theta}}.$$

Also, for $\mathbf{u} \in \mathbb{R}^p$,

$$
\begin{aligned}
(3.14) \quad &\tilde{W}_{m_n}^\dagger(\mathbf{u}) \\
&:= \sum_{t=p+1}^{m_n}\left[\ln f\left(Z_t^* + m_n^{-1/\alpha_0}\sum_{j\neq 0}\hat{c}_j(\mathbf{u})Z_{t-j}^*; \boldsymbol{\tau}_0\right) - \ln f(Z_t^*; \boldsymbol{\tau}_0)\right]
\end{aligned}
$$

and

$$(3.15) \qquad \tilde{W}_{m_n}(\mathbf{u}) := \mathcal{L}_{m_n}^*(\hat{\boldsymbol{\theta}}_{\mathrm{ML}} + m_n^{-1/\alpha_0}\mathbf{u}, s_0) - \mathcal{L}_{m_n}^*(\hat{\boldsymbol{\theta}}_{\mathrm{ML}}, s_0).$$

Now, let $\mathcal{M}_p(C(\mathbb{R}^p))$ represent the space of probability measures on $C(\mathbb{R}^p)$, and let $d_0$ metrize the topology of weak convergence on $\mathcal{M}_p(C(\mathbb{R}^p))$. That is, for random elements $L_n$ and $L$ of $\mathcal{M}_p(C(\mathbb{R}^p))$, $L_n \overset{P}{\to} L$ if and only if $d_0(L_n, L) \overset{P}{\to} 0$ on $\mathbb{R}$, and there exists a dense sequence $\{\tilde{h}_j\}_{j=1}^{\infty}$ of bounded, continuous functions on $C(\mathbb{R}^p)$ such that $d_0(L_n, L) \overset{P}{\to} 0$ is equivalent to $\int_{C(\mathbb{R}^p)} \tilde{h}_j \, dL_n \overset{P}{\to} \int_{C(\mathbb{R}^p)} \tilde{h}_j \, dL$ on $\mathbb{R}$ for all $j \in \{1, 2, \ldots\}$. We now show that, if $L_n(\cdot) := P(\tilde{W}_{m_n} \in \cdot | X_1, \ldots, X_n)$ and $L_n^\dagger(\cdot) := P(\tilde{W}_{m_n}^\dagger \in \cdot | X_1, \ldots, X_n)$, then $L_n - L_n^\dagger \overset{P}{\to} 0$ on $\mathcal{M}_p(C(\mathbb{R}^p))$. Following the proof of Theorem 2.1 in [15], it suffices to show that for any subsequence $\{n_k\}$ there exists a further subsequence $\{n_{k'}\}$ for which $L_{n_{k'}} - L_{n_{k'}}^\dagger \overset{\text{a.s.}}{\to} 0$ relative to the metric $d_0$, which holds if, for almost all realizations of $\{X_t\}$, $\tilde{W}_{m_{n_{k'}}}(\cdot) - \tilde{W}_{m_{n_{k'}}}^\dagger(\cdot) \overset{P}{\to} 0$ on $C(\mathbb{R}^p)$. By Lemma A.13, for any subsequence, any $T \in \{1, 2, \ldots\}$ and any $\kappa \in \{1, 1/2, 1/3, \ldots\}$, there exists a further subsequence $\{n_{k'}^{T,\kappa}\}$ for which $\mathrm{P}(\sup_{\|\mathbf{u}\| \leq T} |\tilde{W}_{m_{n_{k'}^{T,\kappa}}}(\mathbf{u}) - \tilde{W}_{m_{n_{k'}^{T,\kappa}}}^{\dagger}(\mathbf{u})| > \kappa | X_1, \ldots, X_{n_{k'}^{T,\kappa}}) \overset{\text{a.s.}}{\to} 0$. Using a diagonal sequence argument, it follows that there exists a subsequence $\{n_{k'}\}$ of $\{n_k\}$ for which $\mathrm{P}(\sup_{\|\mathbf{u}\| \leq T} |\tilde{W}_{m_{n_{k'}}}(\mathbf{u}) - \tilde{W}_{m_{n_{k'}}}^{\dagger}(\mathbf{u})| > \kappa | X_1, \ldots, X_{n_{k'}}) \to 0$ for almost all $\{X_t\}$ and any $T, \kappa > 0$ and, thus, $\tilde{W}_{m_{n_{k'}}}(\cdot) - \tilde{W}_{m_{n_{k'}}}^{\dagger}(\cdot) \overset{P}{\to} 0$ on $C(\mathbb{R}^p)$ for almost all $\{X_t\}$.

Following the proof of Theorem 3.1 in [15], $L_n^{\dagger}(\cdot) = \mathrm{P}(\tilde{W}_{m_n}^{\dagger} \in \cdot | X_1, \ldots, X_n) \overset{P}{\to} \mathrm{P}(W \in \cdot)$ on $\mathcal{M}_p(C(\mathbb{R}^p))$, and so $L_n(\cdot) = \mathrm{P}(\tilde{W}_{m_n} \in \cdot | X_1, \ldots, X_n) \overset{P}{\to} \mathrm{P}(W \in \cdot)$ on $\mathcal{M}_p(C(\mathbb{R}^p))$ also. Therefore, because $\mathcal{L}_{m_n}^*(\boldsymbol{\theta}, s_0) - \mathcal{L}_{m_n}^*(\hat{\boldsymbol{\theta}}_{\mathrm{ML}}, s_0) = \tilde{W}_{m_n}(m_n^{1/\alpha_0} \times (\boldsymbol{\theta} - \hat{\boldsymbol{\theta}}_{\mathrm{ML}}))$ and $\boldsymbol{\xi}$ uniquely maximizes $W(\cdot)$ almost surely, it can be shown that there exists a sequence of maximizers $\hat{\boldsymbol{\theta}}_{m_n}^*$ of $\mathcal{L}_{m_n}^*(\cdot, s_0)$, such that $\mathrm{P}(m_n^{1/\alpha_0}(\hat{\boldsymbol{\theta}}_{m_n}^* - \hat{\boldsymbol{\theta}}_{\mathrm{ML}}) \in \cdot | X_1 \ldots, X_n) \overset{P}{\to} \mathrm{P}(\boldsymbol{\xi} \in \cdot)$ on $\mathcal{M}_p(\mathbb{R}^p)$ (the proof is similar to that of Theorem 2.2 in [15]). Since

$$m_n^{1/\hat{\alpha}_{\mathrm{ML}}}(\hat{\boldsymbol{\theta}}_{m_n}^* - \hat{\boldsymbol{\theta}}_{\mathrm{ML}}) - m_n^{1/\alpha_0}(\hat{\boldsymbol{\theta}}_{m_n}^* - \hat{\boldsymbol{\theta}}_{\mathrm{ML}})$$
$$= -\left(\frac{m_n^{1/\alpha_n^*} \ln(m_n)}{m_n^{1/\alpha_0}(\alpha_n^*)^2}\right)(\hat{\alpha}_{\mathrm{ML}} - \alpha_0) m_n^{1/\alpha_0}(\hat{\boldsymbol{\theta}}_{m_n}^* - \hat{\boldsymbol{\theta}}_{\mathrm{ML}}),$$

where $\alpha_n^*$ lies between $\hat{\alpha}_{\mathrm{ML}}$ and $\alpha_0$, and $n^{1/2}(\hat{\alpha}_{\mathrm{ML}} - \alpha_0) = O_p(1)$, $\mathrm{P}(\|(m_n^{1/\hat{\alpha}_{\mathrm{ML}}} - m_n^{1/\alpha_0})(\hat{\boldsymbol{\theta}}_{m_n}^* - \hat{\boldsymbol{\theta}}_{\mathrm{ML}})\| > \kappa | X_1, \ldots, X_n) \overset{P}{\to} 0$ for any $\kappa > 0$. Hence, $\mathrm{P}(m_n^{1/\hat{\alpha}_{\mathrm{ML}}}(\hat{\boldsymbol{\theta}}_{m_n}^* - \hat{\boldsymbol{\theta}}_{\mathrm{ML}}) \in \cdot | X_1, \ldots, X_n) \overset{P}{\to} \mathrm{P}(\boldsymbol{\xi} \in \cdot)$ on $\mathcal{M}_p(\mathbb{R}^p)$. The mean-value theorem can be used to show that (3.12) holds. $\square$

Thus, $m_n^{1/\hat{\alpha}_{\mathrm{ML}}}(\hat{\boldsymbol{\theta}}_{m_n}^* - \hat{\boldsymbol{\theta}}_{\mathrm{ML}})$ and $m_n^{1/\hat{\alpha}_{\mathrm{ML}}}(\hat{\boldsymbol{\phi}}_{m_n}^* - \hat{\boldsymbol{\phi}}_{\mathrm{ML}})$, conditioned on $\{X_t\}_{t=1}^n$, have the same limiting distributions as $n^{1/\alpha_0}(\hat{\boldsymbol{\theta}}_{\mathrm{ML}} - \boldsymbol{\theta}_0)$ and $n^{1/\alpha_0}(\hat{\boldsymbol{\phi}}_{\mathrm{ML}} - \boldsymbol{\phi}_0)$, respectively. If $n$ is large, these limiting distributions can, therefore, be approximated by simulating bootstrap values of $\hat{\boldsymbol{\theta}}_{m_n}^*$ and $\hat{\boldsymbol{\phi}}_{m_n}^*$, and looking at the distributions for $m_n^{1/\hat{\alpha}_{\mathrm{ML}}}(\hat{\boldsymbol{\theta}}_{m_n}^* - \hat{\boldsymbol{\theta}}_{\mathrm{ML}})$ and $m_n^{1/\hat{\alpha}_{\mathrm{ML}}}(\hat{\boldsymbol{\phi}}_{m_n}^* - \hat{\boldsymbol{\phi}}_{\mathrm{ML}})$. In principle, one could also examine the limiting distributions for $n^{1/\alpha_0}(\hat{\boldsymbol{\theta}}_{\mathrm{ML}} - \boldsymbol{\theta}_0)$ and $n^{1/\alpha_0}(\hat{\boldsymbol{\phi}}_{\mathrm{ML}} - \boldsymbol{\phi}_0)$ by simulating realizations of $W(\cdot)$, with the true parameter values $\boldsymbol{\theta}_0$ and $\boldsymbol{\tau}_0$ replaced by estimates, and by finding the corresponding values of the maximizer $\boldsymbol{\xi}$, but this procedure is much more laborious than the bootstrap. Confidence intervals for the elements of $\boldsymbol{\theta}_0$ and $\boldsymbol{\phi}_0$ can be obtained using the limiting results for $\hat{\boldsymbol{\theta}}_{\mathrm{ML}}$ and $\hat{\boldsymbol{\phi}}_{\mathrm{ML}}$ in (3.5) and (3.6), bootstrap estimates of quantiles for the limiting distributions and the estimate $\hat{\alpha}_{\mathrm{ML}}$ of $\alpha_0$.

For the elements of $\boldsymbol{\tau}_0$, confidence intervals can be directly obtained from the limiting result for $\hat{\boldsymbol{\tau}}_{\mathrm{ML}}$ in (3.5). Because $\mathbf{I}^{-1}(\cdot)$ is continuous at $\boldsymbol{\tau}_0$ and $\hat{\boldsymbol{\tau}}_{\mathrm{ML}} \xrightarrow{P} \boldsymbol{\tau}_0$, $\mathbf{I}^{-1}(\hat{\boldsymbol{\tau}}_{\mathrm{ML}})$ is a consistent estimator for $\mathbf{I}^{-1}(\boldsymbol{\tau}_0)$ which can be used to compute standard errors for the estimates.

## 4. Numerical results.

4.1. *Simulation study.* In this section we describe a simulation experiment to study the behavior of the ML estimators for finite samples. We did these simulations in MATLAB, using John Nolan's STABLE library (http://academic2.american.edu/~jpnolan/stable/stable.html) to generate stable noise and evaluate stable densities. The STABLE library uses the algorithm in Chambers, Mallows and Stuck [9] to generate stable noise and the algorithm in Nolan [30] to evaluate stable densities.

For each of 300 replicates, we simulated an AR series of length $n = 500$ with stable noise and then found $\hat{\boldsymbol{\eta}}_{\mathrm{ML}} = (\hat{\boldsymbol{\theta}}'_{\mathrm{ML}}, \hat{\boldsymbol{\tau}}'_{\mathrm{ML}})'$ by maximizing the log-likelihood $\mathcal{L}$ in (2.12) with respect to both $s \in \{0, \ldots, p\}$ and $\boldsymbol{\eta}$. To reduce the possibility of the optimizer getting trapped at local maxima, for each $s \in \{0, \ldots, p\}$, we used 1200 randomly chosen starting values for $\boldsymbol{\eta}$. We evaluated the log-likelihood at each of the candidate values and, for each $s \in \{0, \ldots, p\}$, reduced the collection of initial values to the eight with the highest likelihoods. Optimized values were found using the Nelder–Mead algorithm (see, e.g., Lagarias et al. [23]) and the $8(p + 1)$ initial values as starting points. The optimized value for which the likelihood was highest was chosen to be $\hat{\boldsymbol{\eta}}_{\mathrm{ML}}$, and then $\hat{\boldsymbol{\phi}}_{\mathrm{ML}}$ was computed using (2.13). In all cases, $\mathcal{L}$ was maximized at $s = s_0$, so the true order of noncausality for the AR model was always correctly identified.

We obtained simulation results for the causal AR(1) model with parameter $\phi_0 = 0.5$, the noncausal AR(1) model with parameter $\phi_0 = 2.0$ and the AR(2) model with parameter $\phi_0 = (-1.2, 1.6)'$. The AR(2) polynomial $1 + 1.2z - 1.6z^2$ equals $(1 - 0.8z)(1 + 2z)$, and so it has one root inside and the other outside the unit circle. Results of the simulations appear in Table 1, where we give the empirical means and standard deviations for the parameter estimates. The asymptotic standard deviations were obtained using Theorem 3.2 and values for $\mathbf{I}^{-1}(\boldsymbol{\tau}_0)$ in Nolan [31]. (Values for $\mathbf{I}^{-1}(\cdot)$ not given in Nolan [31] can be computed using the STABLE library.) Results for symmetric stable noise are given on the left-hand side of the table, and results for asymmetric stable noise with $\beta_0 = 0.5$ are given on the right-hand side. In Table 1, we see that the MLEs are all approximately unbiased and that the asymptotic standard deviations fairly accurately reflect the true variability of the estimates $\hat{\alpha}_{\mathrm{ML}}$, $\hat{\beta}_{\mathrm{ML}}$, $\hat{\sigma}_{\mathrm{ML}}$, and $\hat{\mu}_{\mathrm{ML}}$. Note that the values of $\hat{\boldsymbol{\phi}}_{\mathrm{ML}}$, $\hat{\alpha}_{\mathrm{ML}}$, $\hat{\beta}_{\mathrm{ML}}$, and $\hat{\mu}_{\mathrm{ML}}$ are less disperse when the noise distribution is

Table 1

*Empirical means and standard deviations for* ML *estimates of* AR *model parameters. The asymptotic standard deviations were computed using Theorem 3.2 and Nolan [31]*

| | Asymp. std. dev. | Empirical mean | Empirical std. dev. | | Asymp. std. dev. | Empirical mean | Empirical std. dev. |
|---|---|---|---|---|---|---|---|
| $\phi_{01} = 0.5$ | | 0.500 | 0.001 | $\phi_{01} = 0.5$ | | 0.500 | 0.001 |
| $\alpha_0 = 0.8$ | 0.051 | 0.795 | 0.040 | $\alpha_0 = 0.8$ | 0.049 | 0.799 | 0.035 |
| $\beta_0 = 0.0$ | 0.067 | 0.000 | 0.064 | $\beta_0 = 0.5$ | 0.058 | 0.504 | 0.060 |
| $\sigma_0 = 1.0$ | 0.077 | 0.996 | 0.068 | $\sigma_0 = 1.0$ | 0.074 | 0.995 | 0.075 |
| $\mu_0 = 0.0$ | 0.054 | 0.003 | 0.057 | $\mu_0 = 0.0$ | 0.062 | $-0.002$ | 0.066 |
| $\phi_{01} = 0.5$ | | 0.498 | 0.019 | $\phi_{01} = 0.5$ | | 0.500 | 0.018 |
| $\alpha_0 = 1.5$ | 0.071 | 1.499 | 0.069 | $\alpha_0 = 1.5$ | 0.070 | 1.500 | 0.066 |
| $\beta_0 = 0.0$ | 0.137 | 0.012 | 0.142 | $\beta_0 = 0.5$ | 0.121 | 0.491 | 0.121 |
| $\sigma_0 = 1.0$ | 0.048 | 0.997 | 0.050 | $\sigma_0 = 1.0$ | 0.047 | 0.996 | 0.047 |
| $\mu_0 = 0.0$ | 0.078 | $-0.002$ | 0.074 | $\mu_0 = 0.0$ | 0.078 | 0.005 | 0.082 |
| $\phi_{01} = 2.0$ | | 2.000 | 0.004 | $\phi_{01} = 2.0$ | | 2.000 | 0.004 |
| $\alpha_0 = 0.8$ | 0.051 | 0.797 | 0.041 | $\alpha_0 = 0.8$ | 0.049 | 0.795 | 0.037 |
| $\beta_0 = 0.0$ | 0.067 | 0.000 | 0.066 | $\beta_0 = 0.5$ | 0.058 | 0.499 | 0.060 |
| $\sigma_0 = 1.0$ | 0.077 | 1.004 | 0.072 | $\sigma_0 = 1.0$ | 0.074 | 0.996 | 0.072 |
| $\mu_0 = 0.0$ | 0.054 | 0.004 | 0.055 | $\mu_0 = 0.0$ | 0.062 | 0.000 | 0.063 |
| $\phi_{01} = 2.0$ | | 2.003 | 0.074 | $\phi_{01} = 2.0$ | | 2.013 | 0.073 |
| $\alpha_0 = 1.5$ | 0.071 | 1.505 | 0.074 | $\alpha_0 = 1.5$ | 0.070 | 1.497 | 0.069 |
| $\beta_0 = 0.0$ | 0.137 | 0.008 | 0.138 | $\beta_0 = 0.5$ | 0.121 | 0.504 | 0.119 |
| $\sigma_0 = 1.0$ | 0.048 | 1.000 | 0.056 | $\sigma_0 = 1.0$ | 0.047 | 0.996 | 0.061 |
| $\mu_0 = 0.0$ | 0.078 | $-0.006$ | 0.077 | $\mu_0 = 0.0$ | 0.078 | 0.004 | 0.079 |
| $\phi_{01} = -1.2$ | | $-1.200$ | 0.004 | $\phi_{01} = -1.2$ | | $-1.200$ | 0.004 |
| $\phi_{02} = 1.6$ | | 1.600 | 0.004 | $\phi_{02} = 1.6$ | | 1.600 | 0.004 |
| $\alpha_0 = 0.8$ | 0.051 | 0.798 | 0.041 | $\alpha_0 = 0.8$ | 0.049 | 0.800 | 0.039 |
| $\beta_0 = 0.0$ | 0.067 | $-0.001$ | 0.068 | $\beta_0 = 0.5$ | 0.058 | 0.502 | 0.056 |
| $\sigma_0 = 1.0$ | 0.077 | 0.997 | 0.073 | $\sigma_0 = 1.0$ | 0.074 | 0.997 | 0.071 |
| $\mu_0 = 0.0$ | 0.054 | $-0.002$ | 0.057 | $\mu_0 = 0.0$ | 0.062 | $-0.004$ | 0.064 |
| $\phi_{01} = -1.2$ | | $-1.212$ | 0.083 | $\phi_{01} = -1.2$ | | $-1.204$ | 0.078 |
| $\phi_{02} = 1.6$ | | 1.605 | 0.065 | $\phi_{02} = 1.6$ | | 1.598 | 0.062 |
| $\alpha_0 = 1.5$ | 0.071 | 1.502 | 0.069 | $\alpha_0 = 1.5$ | 0.070 | 1.499 | 0.071 |
| $\beta_0 = 0.0$ | 0.137 | 0.010 | 0.128 | $\beta_0 = 0.5$ | 0.121 | 0.509 | 0.128 |
| $\sigma_0 = 1.0$ | 0.048 | 0.999 | 0.066 | $\sigma_0 = 1.0$ | 0.047 | 0.997 | 0.056 |
| $\mu_0 = 0.0$ | 0.078 | $-0.006$ | 0.078 | $\mu_0 = 0.0$ | 0.078 | 0.000 | 0.083 |

heavier-tailed (ie., when $\alpha_0 = 0.8$), while the values of $\hat{\sigma}_{\mathrm{ML}}$ are more disperse when the noise distribution has heavier tails. Note also that the finite sample results for $\hat{\boldsymbol{\tau}}_{\mathrm{ML}}$ do not appear particularly affected by the value of $\boldsymbol{\phi}_0$, which is not surprising since $\hat{\boldsymbol{\phi}}_{\mathrm{ML}}$ and $\hat{\boldsymbol{\tau}}_{\mathrm{ML}}$ are asymptotically independent.

Normal qq-plots show that, in all cases, $\hat{\alpha}_{\mathrm{ML}}$, $\hat{\beta}_{\mathrm{ML}}$, $\hat{\sigma}_{\mathrm{ML}}$ and $\hat{\mu}_{\mathrm{ML}}$ have approximately Gaussian distributions. To examine the distribution for $n^{1/\alpha_0}(\hat{\boldsymbol{\phi}}_{\mathrm{ML}} - \boldsymbol{\phi}_0)$, in Figure 1, we give kernel estimates for the density of $n^{1/\alpha_0}(\hat{\phi}_{1,\mathrm{ML}} - \phi_{01})$ when $(\phi_{01}, \alpha_0, \beta_0, \sigma_0, \mu_0)$ is $(0.5, 0.8, 0, 1, 0)$, $(0.5, 0.8, 0.5, 1, 0)$, $(0.5, 1.5, 0, 1, 0)$ and $(0.5, 1.5, 0.5, 1, 0)$. For comparison, we also included normal density functions in Figure 1; the means and variances for the normal densities are the corresponding means and variances for the values of $n^{1/\alpha_0}(\hat{\phi}_{1,\mathrm{ML}} - \phi_{01})$. The distribution of $n^{1/\alpha_0}(\hat{\phi}_{1,\mathrm{ML}} - \phi_1)$ appears more peaked and heavier-tailed than Gaussian, but closer to Gaussian as $\alpha_0$ approaches two. Similar behavior is exhibited by other estimators $\hat{\phi}_{j,\mathrm{ML}}$.

4.2. *Autoregressive modeling.* Figure 2 shows the natural logarithms of the volumes of Wal-Mart stock traded daily on the New York Stock Exchange from December 1, 2003 to December 31, 2004. Sample autocorrelation and partial autocorrelation functions for the series are given in Figure 3. Note that, even if a process has infinite second-order moments, the sample correlations and partial correlations can still be useful for identifying a suitable model for the data (see, e.g., Adler, Feldman and Gallagher [1]). Because the sample partial autocorrelation function is approximately zero after lag two and the data appear "spiky," it is reasonable to try modeling this series $\{X_t\}_{t=1}^{274}$ as an AR(2) process with non-Gaussian stable noise. Additionally, Akaike's information criterion (AIC) is smallest at lag two.

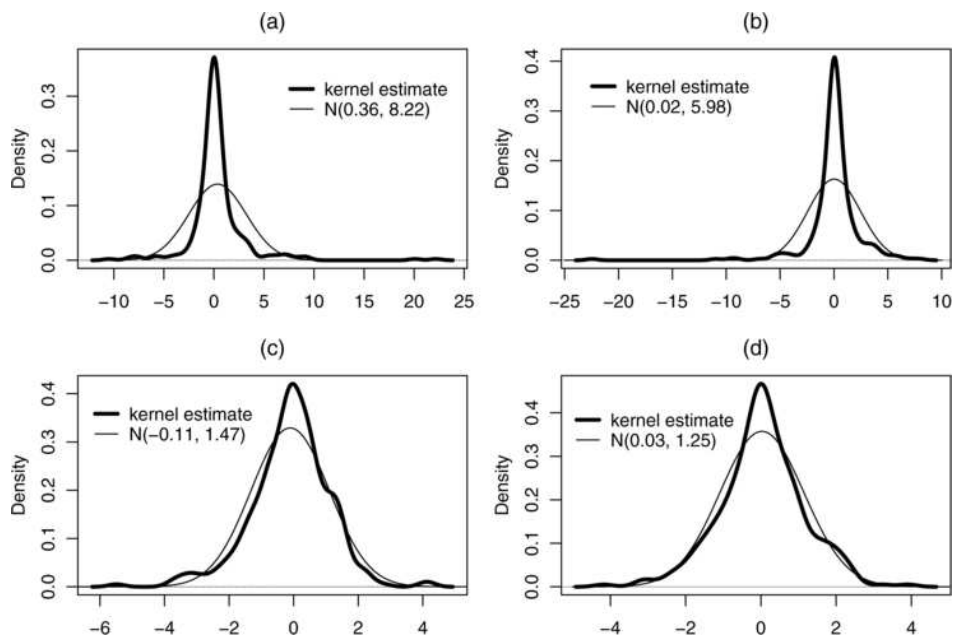

Fig. 1. *Kernel estimates of the density for $n^{1/\alpha_0}(\hat{\phi}_{1,\mathrm{ML}} - \phi_{01})$ when $(\phi_{01}, \alpha_0, \beta_0, \sigma_0, \mu_0)$ is* (a) $(0.5, 0.8, 0, 1, 0)$, (b) $(0.5, 0.8, 0.5, 1, 0)$, (c) $(0.5, 1.5, 0, 1, 0)$ *and* (d) $(0.5, 1.5, 0.5, 1, 0)$, *and normal density functions with the same means and variances as the corresponding values for $n^{1/\alpha_0}(\hat{\phi}_{1,\mathrm{ML}} - \phi_{01})$.*

This supports the suitability of an AR(2) model for $\{X_t\}$. Note that AIC is a consistent order selection criterion for heavy-tailed, infinite variance AR processes (Knight [22]), even though it is not in the finite variance case.

We fit an AR(2) model to $\{X_t\}$ by maximizing $\mathcal{L}$ in (2.12) with respect to both $\boldsymbol{\eta}$ and $s$. The ML estimates are $\hat{\boldsymbol{\eta}}_{\mathrm{ML}} = (\hat{\theta}_1, \hat{\theta}_2, \hat{\alpha}, \hat{\beta}, \hat{\sigma}, \hat{\mu})' = (0.7380, -2.8146, 1.8335, 0.5650, 0.4559, 16.0030)'$, with $s = 1$. Hence, the fit-ted AR(2) polynomial has one root inside and one root outside the unit circle. The residuals from the fitted noncausal AR(2) model

$$(4.1) \quad (1 - 0.7380B)(1 + 2.8146B)X_t = (1 + 2.0766B - 2.0772B^2)X_t = Z_t$$

and sample autocorrelation functions for the absolute values and squares of the mean-corrected residuals are shown in Figure 4(a)–(c). The bounds in (b) and (c) are approximate 95% confidence bounds which we obtained by simulating 100,000 independent sample correlations for the absolute values and squares of 272 mean-corrected i.i.d. stable random variables with $\boldsymbol{\tau} = (1.8335, 0.5650, 0.4559, 16.0030)'$. Based on these graphs, the residuals appear approximately i.i.d., and so we conclude that (4.1) is a satisfactory fitted model for the series $\{X_t\}$. A qq-plot, with empirical quantiles for the residuals plotted against theoretical quantiles of the stable $\boldsymbol{\tau} = (1.8335, 0.5650, 0.4559, 16.0030)'$ distribution, is given in Figure 4(d). Because the qq-plot is remarkably linear, it appears reasonable to model the i.i.d. noise $\{Z_t\}$ in (4.1) as stable with parameter $\boldsymbol{\tau} = (1.8335, 0.5650, 0.4559, 16.0030)'$. Following the discussion at the end of Section 3, approximate 95% bootstrap confidence intervals for $\phi_{01}$ and $\phi_{02}$ are $(-2.2487, -1.8116)$ and $(1.8120, 2.2439)$ (these were obtained from 100 iterations of the bootstrap procedure with $m_n = 135$), and approximate 95% confidence intervals for $\alpha_0$, $\beta_0$, $\sigma_0$ and $\mu_0$, with standard errors computed using $\mathbf{I}^{-1}(\hat{\boldsymbol{\tau}}_{\mathrm{ML}})$, are $(1.6847, 1.9823)$, $(-0.1403, 1)$, $(0.4093, 0.5025)$ and $(15.9102, 16.0958)$.

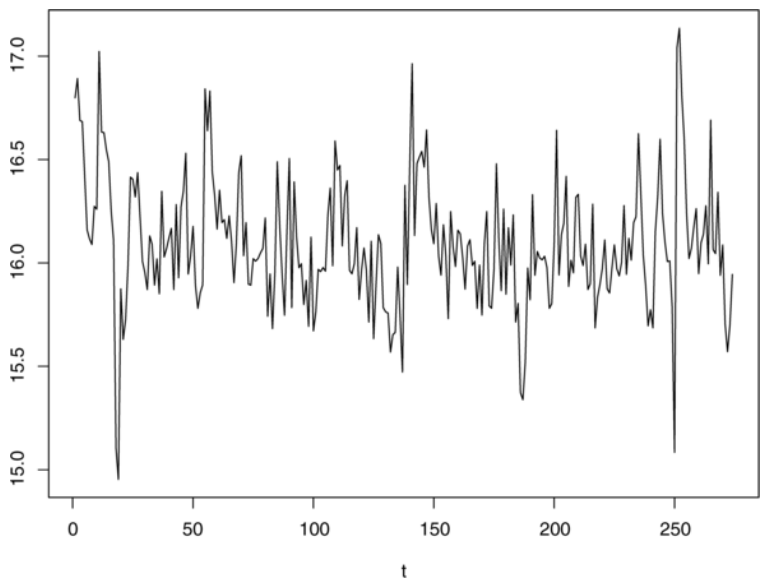

FIG. 2. *The natural logarithms of the volumes of Wal-Mart stock traded daily on the New York Stock Exchange from December 1, 2003 to December 31, 2004.*

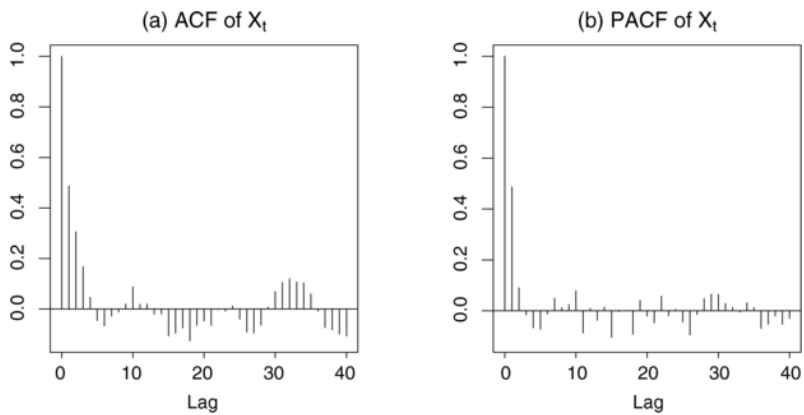

FIG. 3. (a) *The sample autocorrelation function for* $\{X_t\}$ *and* (b) *the sample partial autocorrelation function for* $\{X_t\}$.

In contrast, when we fit a causal AR(2) model to $\{X_t\}$ by maximizing $\mathcal{L}$ with $s = 0$ fixed, we obtain $\hat{\boldsymbol{\eta}} = (\hat{\theta}_1, \hat{\theta}_2, \hat{\alpha}, \hat{\beta}, \hat{\sigma}, \hat{\mu})' = (0.4326, 0.2122, 1.7214, 0.5849, 0.1559, 5.6768)'$. The sample autocorrelation functions for the absolute values and squares of the mean-corrected residuals from this fitted causal model are given in Figure 5. Because both the absolute values and squares have large lag-one correlations, the residuals do not appear independent, and so the causal AR model is not suitable for $\{X_t\}$.

## APPENDIX

In this final section, we give proofs of the lemmas used to establish the results of Section 3.

LEMMA A.1. *For any fixed* $\mathbf{u} \in \mathbb{R}^p$ *and for* $Z_{k,j}(\mathbf{u})$ *between* $Z_{k,j}$ *and* $Z_{k,j} + [\tilde{c}(\alpha_0)]^{1/\alpha_0}\sigma_0 c_j(\mathbf{u})\delta_k \Gamma_k^{-1/\alpha_0}$,

$$(A.1) \quad \sum_{k=1}^{\infty}\sum_{j \neq 0}|c_j(\mathbf{u})|\Gamma_k^{-1/\alpha_0}\left|\frac{\partial \ln f(Z_{k,j}(\mathbf{u}); \boldsymbol{\tau}_0)}{\partial z} - \frac{\partial \ln f(Z_{k,j}; \boldsymbol{\tau}_0)}{\partial z}\right|$$

*is finite a.s.*

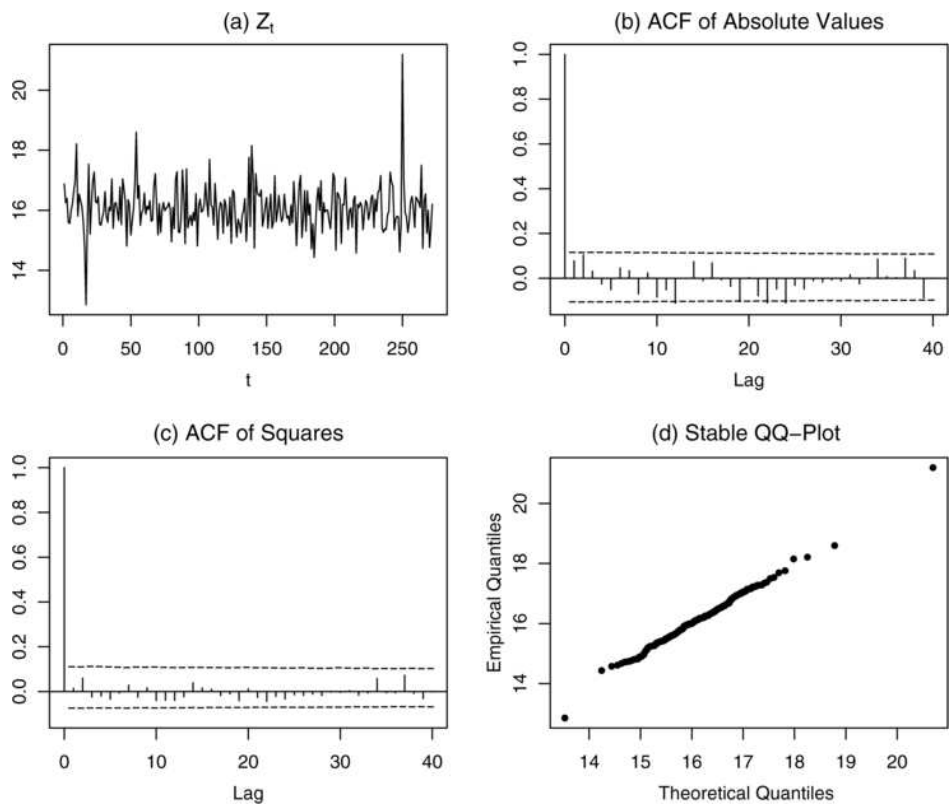

FIG. 4.   (a) *The residuals* $\{Z_t\}$, (b) *the sample autocorrelation function for the absolute values of mean-corrected* $\{Z_t\}$, (c) *the sample autocorrelation function for the squares of mean-corrected* $\{Z_t\}$ *and* (d) *the stable qq-plot for* $\{Z_t\}$.

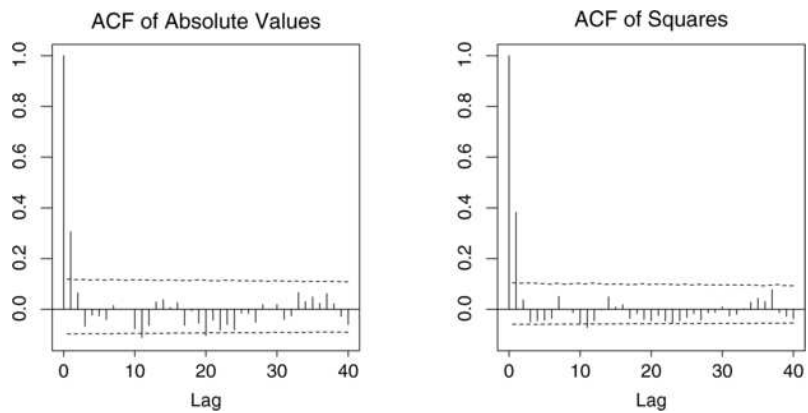

FIG. 5.   *The sample autocorrelation functions for the absolute values and squares of the mean-corrected residuals from the fitted causal AR(2) model.*

PROOF. Since equation (A.1) equals $\sum_{k=1}^{\infty} \sum_{j\neq 0} |c_j(\mathbf{u})| \Gamma_k^{-1/\alpha_0} |Z_{k,j}(\mathbf{u}) - Z_{k,j}| |\partial^2 \ln f(Z_{k,j}^*(\mathbf{u}); \boldsymbol{\tau}_0)/\partial z^2|$, where $Z_{k,j}^*(\mathbf{u})$ is between $Z_{k,j}$ and $Z_{k,j}(\mathbf{u})$, (A.1) is bounded above by

$$(A.2) \qquad [\tilde{c}(\alpha_0)]^{1/\alpha_0} \sigma_0 \sup_{z\in\mathbb{R}} \left| \frac{\partial^2 \ln f(z; \boldsymbol{\tau}_0)}{\partial z^2} \right| \sum_{k=1}^{\infty} \Gamma_k^{-2/\alpha_0} \sum_{j\neq 0} c_j^2(\mathbf{u}).$$

By (2.5) and the continuity of $\partial^2 \ln f(\cdot; \boldsymbol{\tau}_0)/\partial z^2$ on $\mathbb{R}$, $\sup_{z\in\mathbb{R}} |\partial^2 \ln f(z; \boldsymbol{\tau}_0)/\partial z^2| < \infty$. Now suppose $k^\dagger \in \{2, 3, \ldots\}$ such that $k^\dagger > 2/\alpha_0$. It follows that $E\{\sum_{k=k^\dagger}^{\infty} \Gamma_k^{-2/\alpha_0}\} = \sum_{k=k^\dagger}^{\infty} \Gamma(k-2/\alpha_0)/\Gamma(k) < (constant) \sum_{k=k^\dagger}^{\infty} k^{-2/\alpha_0} < \infty$. Consequently, since $0 < [\tilde{c}(\alpha_0)]^{1/\alpha_0} \sigma_0 < \infty$, $\sum_{j\neq 0} c_j^2(\mathbf{u}) < \infty$ by (3.3) and $\sum_{k=1}^{k^\dagger - 1} \Gamma_k^{-2/\alpha_0} < \infty$ a.s., (A.2) is finite a.s. □

LEMMA A.2. *For any fixed* $\mathbf{u} \in \mathbb{R}^p$,

$$(A.3) \qquad \sum_{k=1}^{\infty} \sum_{j\neq 0} \left| c_j(\mathbf{u})(\Gamma_k^{-1/\alpha_0} - k^{-1/\alpha_0}) \frac{\partial \ln f(Z_{k,j}; \boldsymbol{\tau}_0)}{\partial z} \right| < \infty \qquad a.s.$$

PROOF. The left-hand side of (A.3) is bounded above by $\sup_{z\in\mathbb{R}} |\partial \ln f(z; \boldsymbol{\tau}_0)/\partial z| \sum_{k=1}^{\infty} |\Gamma_k^{-1/\alpha_0} - k^{-1/\alpha_0}| \sum_{j\neq 0} |c_j(\mathbf{u})|$. By (2.6), $\sup_{z\in\mathbb{R}} |\partial \ln f(z; \boldsymbol{\tau}_0)/\partial z| < \infty$, by (3.3), $\sum_{j\neq 0} |c_j(\mathbf{u})| < \infty$, and, from the proof of Proposition A.3 in Davis, Knight and Liu [12], $\sum_{k=1}^{\infty} |\Gamma_k^{-1/\alpha_0} - k^{-1/\alpha_0}| < \infty$ a.s. Thus, (A.3) holds. □

LEMMA A.3. *For any fixed* $\mathbf{u} \in \mathbb{R}^p$, $|\sum_{k=1}^{\infty} \sum_{j\neq 0} c_j(\mathbf{u}) \delta_k k^{-1/\alpha_0} [\partial \ln f(Z_{k,j}; \boldsymbol{\tau}_0)/\partial z]| < \infty$ a.s.

PROOF. The sequence $\{\sum_{j\neq 0} c_j(\mathbf{u}) \delta_k k^{-1/\alpha_0} [\partial \ln f(Z_{k,j}; \boldsymbol{\tau}_0)/\partial z]\}_{k=1}^{\infty}$ is a series of independent random variables which, by dominated convergence, all have mean zero, since $\sum_{j\neq 0} |c_j(\mathbf{u})| < \infty$, $\sup_{z\in\mathbb{R}} |\partial \ln f(z; \boldsymbol{\tau}_0)/\partial z| < \infty$ and $E\{\partial \ln f(Z_{k,j}; \boldsymbol{\tau}_0)/\partial z\} = \int_{-\infty}^{\infty} (\partial f(z; \boldsymbol{\tau}_0)/\partial z) \, dz = 0$. Therefore, because

$$\sum_{k=1}^{\infty} \text{Var}\left\{ \sum_{j\neq 0} c_j(\mathbf{u}) \delta_k k^{-1/\alpha_0} \frac{\partial \ln f(Z_{k,j}; \boldsymbol{\tau}_0)}{\partial z} \right\}$$

$$\leq \left( \sup_{z\in\mathbb{R}} \left| \frac{\partial \ln f(z; \boldsymbol{\tau}_0)}{\partial z} \right| \right)^2 \left( \sum_{j\neq 0} |c_j(\mathbf{u})| \right)^2 \sum_{k=1}^{\infty} k^{-2/\alpha_0}$$

$$< \infty,$$

the result holds by the Kolmogorov convergence theorem (see, e.g., Resnick [33], page 212). □

LEMMA A.4. *For $\mathbf{u} \in \mathbb{R}^p$ and $\mathbf{v} \in \mathbb{R}^4$,*

$$
\begin{aligned}
\text{(A.4)} \qquad & \sum_{t=p+1}^{n} \ln f\left(Z_t\left(\boldsymbol{\theta}_0 + \frac{\mathbf{u}}{n^{1/\alpha_0}}, s_0\right); \boldsymbol{\tau}_0 + \frac{\mathbf{v}}{\sqrt{n}}\right) \\
& - \sum_{t=p+1}^{n} \ln f\left(Z_t + n^{-1/\alpha_0} \sum_{j=-\infty}^{\infty} c_j(\mathbf{u}) Z_{t-j}; \boldsymbol{\tau}_0 + \frac{\mathbf{v}}{\sqrt{n}}\right),
\end{aligned}
$$

*with $Z_t(\cdot, \cdot)$ as defined in (2.11), converges in probability to zero on $C(\mathbb{R}^{p+4})$ as $n \to \infty$.*

PROOF. Let $T > 0$. We begin by showing that (A.4) is $o_p(1)$ on $C([-T, T]^{p+4})$. Since $\{Z_t(\boldsymbol{\theta}_0, s_0)\} = \{Z_t\}$, and following (3.2), equation (A.4) equals

$$
\begin{aligned}
\text{(A.5)} \qquad & \sum_{t=p+1}^{n} \left\{ \frac{\partial \ln f(Z_{t,n}^*(\mathbf{u}); \boldsymbol{\tau}_0 + \mathbf{v}/\sqrt{n})}{\partial z} \right. \\
& \left. \times \left[ Z_t\left(\boldsymbol{\theta}_0 + \frac{\mathbf{u}}{n^{1/\alpha_0}}, s_0\right) - Z_t(\boldsymbol{\theta}_0, s_0) - \frac{\mathbf{u}'}{n^{1/\alpha_0}} \frac{\partial Z_t(\boldsymbol{\theta}_0, s_0)}{\partial \boldsymbol{\theta}} \right] \right\},
\end{aligned}
$$

where $Z_{t,n}^*(\mathbf{u})$ lies between $Z_t(\boldsymbol{\theta}_0 + n^{-1/\alpha_0}\mathbf{u}, s_0)$ and $Z_t + n^{-1/\alpha_0}\mathbf{u}'\partial Z_t(\boldsymbol{\theta}_0, s_0)/\partial\boldsymbol{\theta}$. Equation (A.5) can be expressed as

$$
\frac{1}{2n^{2/\alpha_0}} \sum_{t=p+1}^{n} \frac{\partial \ln f(Z_{t,n}^*(\mathbf{u}); \boldsymbol{\tau}_0 + \mathbf{v}/\sqrt{n})}{\partial z} \mathbf{u}' \frac{\partial^2 Z_t(\boldsymbol{\theta}_{t,n}^*(\mathbf{u}), s_0)}{\partial \boldsymbol{\theta} \, \partial \boldsymbol{\theta}'} \mathbf{u},
$$

with $\boldsymbol{\theta}_{t,n}^*(\mathbf{u})$ between $\boldsymbol{\theta}_0$ and $\boldsymbol{\theta}_0 + n^{-1/\alpha_0}\mathbf{u}$. Following (3.1), the mixed partial derivatives of $Z_t(\boldsymbol{\theta}, s)$ are given by

$$
\frac{\partial^2 Z_t(\boldsymbol{\theta}, s)}{\partial \theta_j \, \partial \theta_k} = \begin{cases} 0, & j, k = 1, \ldots, r, \\ X_{t+r-j-k}, & j = 1, \ldots, r, \ k = r+1, \ldots, p, \\ 0, & j, k = r+1, \ldots, p, \end{cases}
$$

and so we have

$$
\begin{aligned}
& \sup_{(\mathbf{u}', \mathbf{v}')' \in [-T,T]^{p+4}} \frac{1}{2n^{2/\alpha_0}} \left| \sum_{t=p+1}^{n} \frac{\partial \ln f(Z_{t,n}^*(\mathbf{u}); \boldsymbol{\tau}_0 + \mathbf{v}/\sqrt{n})}{\partial z} \right. \\
& \qquad\qquad\qquad\qquad \left. \times \mathbf{u}' \frac{\partial^2 Z_t(\boldsymbol{\theta}_{t,n}^*(\mathbf{u}), s_0)}{\partial \boldsymbol{\theta} \, \partial \boldsymbol{\theta}'} \mathbf{u} \right| \\
& \leq \sup_{z \in \mathbb{R}, \mathbf{v} \in [-T,T]^4} \left| \frac{\partial \ln f(z; \boldsymbol{\tau}_0 + \mathbf{v}/\sqrt{n})}{\partial z} \right| \\
& \qquad \times \sup_{\mathbf{u} \in [-T,T]^p} \frac{1}{2n^{2/\alpha_0}} \sum_{t=p+1}^{n} \left| \mathbf{u}' \frac{\partial^2 Z_t(\boldsymbol{\theta}_{t,n}^*(\mathbf{u}), s_0)}{\partial \boldsymbol{\theta} \, \partial \boldsymbol{\theta}'} \mathbf{u} \right|
\end{aligned}
$$

$$\leq \sup_{z\in\mathbb{R},\mathbf{v}\in[-T,T]^4}\left|\frac{\partial\ln f(z;\boldsymbol{\tau}_0+\mathbf{v}/\sqrt{n})}{\partial z}\right|\frac{T^2p^2}{n^{2/\alpha_0}}\sum_{t=p+1}^{n}\sum_{j=2}^{p}|X_{t-j}|$$

$$\leq \sup_{z\in\mathbb{R},\mathbf{v}\in[-T,T]^4}\left|\frac{\partial\ln f(z;\boldsymbol{\tau}_0+\mathbf{v}/\sqrt{n})}{\partial z}\right|\frac{T^2p^2}{n^{2/\alpha_0}}$$

(A.6)

$$\times\sum_{t=p+1}^{n}\sum_{j=2}^{p}\sum_{k=-\infty}^{\infty}|\psi_k Z_{t-j-k}|$$

(recall that $X_t=\sum_{j=-\infty}^{\infty}\psi_j Z_{t-j}$). By (2.6), $\sup_{z\in\mathbb{R},\mathbf{v}\in[-T,T]^4}|\partial\ln f(z;\boldsymbol{\tau}_0+\mathbf{v}/\sqrt{n})/\partial z|=O(1)$ as $n\to\infty$. Now let $\epsilon>0$ and $\kappa_1=(3/4)\alpha_0 I\{\alpha_0\leq 1\}+I\{\alpha_0>1\}$, and observe that $\mathrm{E}|Z_1|^{\kappa_1}<\infty$ and $0<\kappa_1\leq 1$. Using the Markov inequality,

$$\mathrm{P}\Bigg(\Bigg[\frac{1}{n^{2/\alpha_0}}\sum_{t=p+1}^{n}\sum_{j=2}^{p}\sum_{k=-\infty}^{\infty}|\psi_k Z_{t-j-k}|\Bigg]^{\kappa_1}>\epsilon^{\kappa_1}\Bigg)$$

$$\leq\left(\frac{1}{\epsilon n^{2/\alpha_0}}\right)^{\kappa_1}\mathrm{E}\Bigg\{\sum_{t=p+1}^{n}\sum_{j=2}^{p}\sum_{k=-\infty}^{\infty}|\psi_k Z_{t-j-k}|\Bigg\}^{\kappa_1}$$

$$\leq\left(\frac{1}{\epsilon n^{2/\alpha_0}}\right)^{\kappa_1}\mathrm{E}\Bigg\{\sum_{t=p+1}^{n}\sum_{j=2}^{p}\sum_{k=-\infty}^{\infty}|\psi_k Z_{t-j-k}|^{\kappa_1}\Bigg\}$$

$$\leq\epsilon^{-\kappa_1}n^{1-2\kappa_1/\alpha_0}p\mathrm{E}|Z_1|^{\kappa_1}\sum_{k=-\infty}^{\infty}|\psi_k|^{\kappa_1}$$

$$\overset{n\to\infty}{\to}0.$$

Consequently, (A.6) is $o_p(1)$ on $\mathbb{R}$, and so (A.4) is $o_p(1)$ on $C([-T,T]^{p+4})$. Since $T>0$ was arbitrarily chosen, for any compact set $K\subset\mathbb{R}^{p+4}$, (A.4) is $o_p(1)$ on $C(K)$, and it therefore follows that (A.4) is $o_p(1)$ on $C(\mathbb{R}^{p+4})$. $\square$

LEMMA A.5. *For* $\mathbf{u}\in\mathbb{R}^p$ *and* $\mathbf{v}\in\mathbb{R}^4$,

$$\sum_{t=p+1}^{n}\ln f\Bigg(Z_t+n^{-1/\alpha_0}\sum_{j=-\infty}^{\infty}c_j(\mathbf{u})Z_{t-j};\boldsymbol{\tau}_0+\frac{\mathbf{v}}{\sqrt{n}}\Bigg)$$

(A.7)

$$-\sum_{t=p+1}^{n}\ln f\Bigg(Z_t+n^{-1/\alpha_0}\sum_{j=-\infty}^{\infty}c_j(\mathbf{u})Z_{t-j};\boldsymbol{\tau}_0\Bigg)$$

$$-\frac{\mathbf{v}'}{\sqrt{n}}\sum_{t=p+1}^{n}\frac{\partial\ln f(Z_t;\boldsymbol{\tau}_0)}{\partial\boldsymbol{\tau}}+\frac{1}{2}\mathbf{v}'\mathbf{I}(\boldsymbol{\tau}_0)\mathbf{v}$$

*converges in probability to zero on* $C(\mathbb{R}^{p+4})$ *as* $n\to\infty$.

Proof. Using a Taylor series expansion about $\boldsymbol{\tau}_0$, equation (A.7) equals

$$
\text{(A.8)} \quad \frac{\mathbf{v}'}{\sqrt{n}} \sum_{t=p+1}^{n} \Bigg[ \frac{\partial \ln f(Z_t + n^{-1/\alpha_0} \sum_{j=-\infty}^{\infty} c_j(\mathbf{u}) Z_{t-j}; \boldsymbol{\tau}_0)}{\partial \boldsymbol{\tau}}
$$

$$
- \frac{\partial \ln f(Z_t; \boldsymbol{\tau}_0)}{\partial \boldsymbol{\tau}} \Bigg]
$$

$$
\text{(A.9)} \quad + \frac{\mathbf{v}'}{2n} \sum_{t=p+1}^{n} \frac{\partial^2 \ln f(Z_t + n^{-1/\alpha_0} \sum_{j=-\infty}^{\infty} c_j(\mathbf{u}) Z_{t-j}; \boldsymbol{\tau}_n^*(\mathbf{v}))}{\partial \boldsymbol{\tau} \, \partial \boldsymbol{\tau}'} \mathbf{v}
$$

$$
+ \frac{1}{2} \mathbf{v}' \mathbf{I}(\boldsymbol{\tau}_0) \mathbf{v},
$$

where $\boldsymbol{\tau}_n^*(\mathbf{v})$ is between $\boldsymbol{\tau}_0$ and $\boldsymbol{\tau}_0 + \mathbf{v}/\sqrt{n}$. Let $T > 0$. We will show that $\sup_{\mathbf{u} \in [-T,T]^p}$ of

$$
\text{(A.10)} \quad \Bigg| \frac{1}{\sqrt{n}} \sum_{t=p+1}^{n} \Bigg[ \frac{\partial \ln f(Z_t + n^{-1/\alpha_0} \sum_{j=-\infty}^{\infty} c_j(\mathbf{u}) Z_{t-j}; \boldsymbol{\tau}_0)}{\partial \alpha}
$$

$$
- \frac{\partial \ln f(Z_t; \boldsymbol{\tau}_0)}{\partial \alpha} \Bigg] \Bigg|
$$

is $o_p(1)$. It can be shown similarly that $\sup_{(\mathbf{u}', \mathbf{v}')' \in [-T,T]^{p+4}}$ of (A.8) is $o_p(1)$, and, using the ergodic theorem, $\sup_{(\mathbf{u}', \mathbf{v}')' \in [-T,T]^{p+4}}$ of (A.9) is $o_p(1)$. Since $T > 0$ was arbitrarily chosen, it follows that (A.7) is $o_p(1)$ on $C(\mathbb{R}^{p+4})$.

Observe that $\sup_{\mathbf{u} \in [-T,T]^p}$ of (A.10) equals

$$
\text{(A.11)} \quad \sup_{\mathbf{u} \in [-T,T]^p} \Bigg| \frac{1}{n^{1/2+1/\alpha_0}} \sum_{t=p+1}^{n} \frac{\partial^2 \ln f(Z_{t,n}^{**}(\mathbf{u}); \boldsymbol{\tau}_0)}{\partial z \, \partial \alpha} \sum_{j=-\infty}^{\infty} c_j(\mathbf{u}) Z_{t-j} \Bigg|,
$$

where $Z_{t,n}^{**}(\mathbf{u})$ is between $Z_t$ and $Z_t + n^{-1/\alpha_0} \sum_{j=-\infty}^{\infty} c_j(\mathbf{u}) Z_{t-j}$. Following (3.2), there must exist constants $C_1 > 0$ and $0 < D_1 < 1$ such that

$$
\text{(A.12)} \quad \sup_{\mathbf{u} \in [-T,T]^p} |c_j(\mathbf{u})| \leq C_1 D_1^{|j|} \qquad \forall j \in \{\dots, -1, 0, 1, \dots\},
$$

and so (A.11) is bounded above by

$$
\text{(A.13)} \quad \sup_{z \in \mathbb{R}} \Bigg| \frac{\partial^2 \ln f(z; \boldsymbol{\tau}_0)}{\partial z \, \partial \alpha} \Bigg| \frac{C_1}{n^{1/2+1/\alpha_0}} \sum_{t=p+1}^{n} \sum_{j=-\infty}^{\infty} D_1^{|j|} |Z_{t-j}|.
$$

By (2.7), $\sup_{z \in \mathbb{R}} |\partial^2 \ln f(z; \boldsymbol{\tau}_0)/(\partial z \, \partial \alpha)| < \infty$. Now let $\epsilon > 0$ and $\kappa_2 = \alpha_0(1 + \alpha_0/3)/(1 + \alpha_0/2)I\{\alpha_0 \leq 1\} + I\{\alpha_0 > 1\}$, so that $\kappa_2(1/2 + 1/\alpha_0) > 1$, $\mathrm{E}|Z_1|^{\kappa_2} <$

$\infty$ and $0 < \kappa_2 \le 1$. Since

$$P\left(\left[\frac{1}{n^{1/2+1/\alpha_0}} \sum_{t=p+1}^{n} \sum_{j=-\infty}^{\infty} D_1^{|j|} |Z_{t-j}|\right]^{\kappa_2} > \epsilon^{\kappa_2}\right)$$

$$\le \epsilon^{-\kappa_2} n^{1-\kappa_2(1/2+1/\alpha_0)} \mathrm{E}|Z_1|^{\kappa_2} \sum_{j=-\infty}^{\infty} (D_1^{\kappa_2})^{|j|}$$

and $n^{1-\kappa_2(1/2+1/\alpha_0)} \to 0$, (A.13) is $o_p(1)$ and therefore $\sup_{\mathbf{u}\in[-T,T]^p}$ of (A.10) must also be $o_p(1)$.  □

LEMMA A.6.  *For* $\mathbf{u} \in \mathbb{R}^p$,

$$
\begin{aligned}
(A.14) \quad & \sum_{t=p+1}^{n} \ln f\left(Z_t + n^{-1/\alpha_0} \sum_{j=-\infty}^{\infty} c_j(\mathbf{u}) Z_{t-j}; \boldsymbol{\tau}_0\right) \\
& - \sum_{t=p+1}^{n} \ln f\left(Z_t + n^{-1/\alpha_0} \sum_{j\neq 0} c_j(\mathbf{u}) Z_{t-j}; \boldsymbol{\tau}_0\right) \\
& - \sum_{t=p+1}^{n} \left[\ln f\left(Z_t + \frac{c_0(\mathbf{u})}{n^{1/\alpha_0}} Z_t; \boldsymbol{\tau}_0\right) - \ln f(Z_t; \boldsymbol{\tau}_0)\right]
\end{aligned}
$$

*converges in probability to zero on* $C(\mathbb{R}^p)$ *as* $n \to \infty$.

PROOF.  Equation (A.14) equals

$$(A.15) \quad \frac{c_0(\mathbf{u})}{n^{1/\alpha_0}} \sum_{t=p+1}^{n} Z_t\left[\frac{\partial \ln f(\tilde{Z}^*_{t,n}(\mathbf{u}); \boldsymbol{\tau}_0)}{\partial z} - \frac{\partial \ln f(\tilde{Z}^{**}_{t,n}(\mathbf{u}); \boldsymbol{\tau}_0)}{\partial z}\right],$$

where $\tilde{Z}^*_{t,n}(\mathbf{u})$ is between $Z_t + n^{-1/\alpha_0} \sum_{j=-\infty}^{\infty} c_j(\mathbf{u}) Z_{t-j}$ and $Z_t + n^{-1/\alpha_0} \times \sum_{j\neq 0} c_j(\mathbf{u}) Z_{t-j}$, and $\tilde{Z}^{**}_{t,n}(\mathbf{u})$ is between $Z_t$ and $(1 + n^{-1/\alpha_0} c_0(\mathbf{u})) Z_t$. For $T > 0$, $\sup_{\mathbf{u}\in[-T,T]^p}$ of the absolute value of (A.15) is bounded above by

$$(A.16) \quad \sup_{\mathbf{u}\in[-T,T]^p} \left|\frac{c_0(\mathbf{u})}{n^{2/\alpha_0}} \sum_{t=p+1}^{n} Z_t \sum_{j\neq 0} c_j(\mathbf{u}) Z_{t-j} \frac{\partial^2 \ln f(Z^{***}_{t,n}(\mathbf{u}); \boldsymbol{\tau}_0)}{\partial z^2}\right|$$

$$(A.17) \quad + \sup_{\mathbf{u}\in[-T,T]^p} \frac{c_0^2(\mathbf{u})}{n^{2/\alpha_0}} \sum_{t=p+1}^{n} \left|Z_t^2 \frac{\partial^2 \ln f(Z^{***}_{t,n}(\mathbf{u}); \boldsymbol{\tau}_0)}{\partial z^2}\right|,$$

where $Z^{***}_{t,n}(\mathbf{u}) = Z_t + n^{-1/\alpha_0} \lambda^\dagger_{t,n}(\mathbf{u}) c_0(\mathbf{u}) Z_t + n^{-1/\alpha_0} \lambda^*_{t,n}(\mathbf{u}) \sum_{j\neq 0} c_j(\mathbf{u}) Z_{t-j}$ for some $\lambda^\dagger_{t,n}(\mathbf{u}), \lambda^*_{t,n}(\mathbf{u}) \in [0,1]$. To complete the proof, we show that (A.16) and (A.17) are $o_p(1)$.

Following (A.12), equation (A.16) is bounded above by $\sup_{z \in \mathbb{R}} |\partial^2 \ln f(z; \boldsymbol{\tau}_0) / \partial z^2| n^{-2/\alpha_0} C_1^2 \sum_{t=p+1}^{n} |Z_t| \sum_{j \neq 0} D_1^{|j|} |Z_{t-j}|$. If $\kappa_3 := (3/4)\alpha_0 I\{\alpha_0 \leq 1\} + I\{\alpha_0 > 1\}$, then, for any $\epsilon > 0$,

$$\mathrm{P}\left( \left[ \frac{1}{n^{2/\alpha_0}} \sum_{t=p+1}^{n} |Z_t| \sum_{j \neq 0} D_1^{|j|} |Z_{t-j}| \right]^{\kappa_3} > \epsilon^{\kappa_3} \right)$$

$$\leq \epsilon^{-\kappa_3} n^{1-2\kappa_3/\alpha_0} (\mathrm{E}|Z_1|^{\kappa_3})^2 \sum_{j \neq 0} (D_1^{\kappa_3})^{|j|},$$

which is $o(1)$, and thus (A.16) is $o_p(1)$.

Equation (A.17) is bounded above by

$$(\text{A.18}) \quad \sup_{\mathbf{u} \in [-T,T]^p} \frac{c_0^2(\mathbf{u})}{n^{2/\alpha_0}} \sum_{t=p+1}^{n} Z_t^2 \left| \frac{\partial^2 \ln f([1 + \lambda_{t,n}^\dagger(\mathbf{u}) c_0(\mathbf{u}) / n^{1/\alpha_0}] Z_t; \boldsymbol{\tau}_0)}{\partial z^2} \right|$$

$$+ \sup_{\mathbf{u} \in [-T,T]^p} \frac{c_0^2(\mathbf{u})}{n^{2/\alpha_0}}$$

$$(\text{A.19}) \qquad\qquad \times \sum_{t=p+1}^{n} Z_t^2 \left| \frac{\partial^2 \ln f(Z_{t,n}^{***}(\mathbf{u}); \boldsymbol{\tau}_0)}{\partial z^2} \right.$$

$$\left. - \frac{\partial^2 \ln f([1 + \lambda_{t,n}^\dagger(\mathbf{u}) c_0(\mathbf{u}) / n^{1/\alpha_0}] Z_t; \boldsymbol{\tau}_0)}{\partial z^2} \right|,$$

and (A.18) is bounded above by $\sup_{z \in \mathbb{R}} |z^2 [\partial^2 \ln f(z; \boldsymbol{\tau}_0) / \partial z^2]| \times \sup_{\mathbf{u} \in [-T,T]^p} n^{-2/\alpha_0} c_0^2(\mathbf{u}) \sum_{t=p+1}^{n} (1 + n^{-1/\alpha_0} \lambda_{t,n}^\dagger(\mathbf{u}) c_0(\mathbf{u}))^{-2}$. Since $n^{1-2/\alpha_0} \to 0$, $\sup_{\mathbf{u} \in [-T,T]^p} |c_0(\mathbf{u})| < \infty$ and, from (2.5), $\sup_{z \in \mathbb{R}} |z^2 [\partial^2 \ln f(z; \boldsymbol{\tau}_0) / \partial z^2]| < \infty$, (A.18) is $o_p(1)$. An upper bound for (A.19) is

$$\sup_{\mathbf{u} \in [-T,T]^p} \frac{c_0^2(\mathbf{u})}{n^{3/\alpha_0}} \sum_{t=p+1}^{n} Z_t^2 \left| \frac{\partial^3 \ln f(\tilde{Z}_{t,n}(\mathbf{u}); \boldsymbol{\tau}_0)}{\partial z^3} \sum_{j \neq 0} c_j(\mathbf{u}) Z_{t-j} \right|$$

$$\leq \sup_{z \in \mathbb{R}} \left| \frac{\partial^3 \ln f(z; \boldsymbol{\tau}_0)}{\partial z^3} \right| \left( \frac{C_1}{n^{1/\alpha_0}} \right)^3 \sum_{t=p+1}^{n} Z_t^2 \sum_{j \neq 0} D_1^{|j|} |Z_{t-j}|,$$

where $\tilde{Z}_{t,n}(\mathbf{u})$ is between $Z_{t,n}^{***}(\mathbf{u})$ and $[1 + \lambda_{t,n}^\dagger(\mathbf{u}) c_0(\mathbf{u}) / n^{1/\alpha_0}] Z_t$. If $\kappa_4 := 3\alpha_0/8$, then, for any $\epsilon > 0$,

$$\mathrm{P}\left( \left[ \frac{1}{n^{3/\alpha_0}} \sum_{t=p+1}^{n} Z_t^2 \sum_{j \neq 0} D_1^{|j|} |Z_{t-j}| \right]^{\kappa_4} > \epsilon^{\kappa_4} \right)$$

$$\leq \epsilon^{-\kappa_4} n^{1-3\kappa_4/\alpha_0} \mathrm{E}\{Z_1^{2\kappa_4}\} \mathrm{E}|Z_1|^{\kappa_4} \sum_{j \neq 0} (D_1^{\kappa_4})^{|j|} \xrightarrow{n \to \infty} 0.$$

Since $\sup_{z \in \mathbb{R}} |\partial^3 \ln f(z; \boldsymbol{\tau}_0)/\partial z^3| < \infty$ (see DuMouchel [17]), it follows that (A.19) is also $o_p(1)$. $\quad\square$

LEMMA A.7.   *For* $\mathbf{u} = (u_1, \dots, u_p)' \in \mathbb{R}^p$,

$$
\begin{aligned}
\text{(A.20)} \quad & \sum_{t=p+1}^{n} \left[ \ln f\left( Z_t + \frac{c_0(\mathbf{u})}{n^{1/\alpha_0}} Z_t; \boldsymbol{\tau}_0 \right) - \ln f(Z_t; \boldsymbol{\tau}_0) \right] \\
& \qquad + (n-p) \ln \left| \frac{\theta_{0p} + n^{-1/\alpha_0} u_p}{\theta_{0p}} \right| I\{s_0 > 0\}
\end{aligned}
$$

*converges in probability to zero on* $C(\mathbb{R}^p)$ *as* $n \to \infty$.

PROOF.   If $s_0 = 0$, the result is trivial since, from (3.2), $c_0(\mathbf{u}) = u_p \theta_{0p}^{-1} I\{s_0 > 0\}$, and so, when $s_0 = 0$, equation (A.20) equals zero for all $\mathbf{u} \in \mathbb{R}^p$. Now consider the case $s_0 > 0$. Choose arbitrary $T > 0$ and note that $\sup_{\mathbf{u} \in [-T,T]^p}$ of the absolute value of (A.20) equals

$$
\begin{aligned}
\text{(A.21)} \quad \sup_{\mathbf{u} \in [-T,T]^p} \Bigg| & \sum_{t=p+1}^{n} \left[ \frac{c_0(\mathbf{u})}{n^{1/\alpha_0}} Z_t \frac{\partial \ln f(Z_t; \boldsymbol{\tau}_0)}{\partial z} \right. \\
& \qquad \left. + \frac{c_0^2(\mathbf{u})}{2 n^{2/\alpha_0}} Z_t^2 \frac{\partial^2 \ln f(Z_{t,n}^\dagger(\mathbf{u}); \boldsymbol{\tau}_0)}{\partial z^2} \right] \\
& \qquad + (n-p) \ln \left| \frac{\theta_{0p} + n^{-1/\alpha_0} u_p}{\theta_{0p}} \right| \Bigg|,
\end{aligned}
$$

where $Z_{t,n}^\dagger(\mathbf{u})$ is between $Z_t$ and $[1 + n^{-1/\alpha_0} c_0(\mathbf{u})] Z_t$. Equation (A.21) is bounded above by

$$
\text{(A.22)} \quad \sup_{\mathbf{u} \in [-T,T]^p} \left| \frac{c_0(\mathbf{u})}{n^{1/\alpha_0}} \sum_{t=p+1}^{n} \left[ 1 + Z_t \frac{\partial \ln f(Z_t; \boldsymbol{\tau}_0)}{\partial z} \right] \right|
$$

$$
\text{(A.23)} \quad + \sup_{\mathbf{u} \in [-T,T]^p} \left| (n-p) \left[ \frac{c_0(\mathbf{u})}{n^{1/\alpha_0}} - \ln \left| \frac{\theta_{0p} + n^{-1/\alpha_0} u_p}{\theta_{0p}} \right| \right] \right|
$$

$$
\text{(A.24)} \quad + \sup_{\mathbf{u} \in [-T,T]^p} \left| \frac{c_0^2(\mathbf{u})}{2 n^{2/\alpha_0}} \sum_{t=p+1}^{n} Z_t^2 \frac{\partial^2 \ln f(Z_{t,n}^\dagger(\mathbf{u}); \boldsymbol{\tau}_0)}{\partial z^2} \right|;
$$

we complete the proof by showing that each of these three terms is $o_p(1)$.

Since $\{1 + Z_t[\partial \ln f(Z_t; \boldsymbol{\tau}_0)/\partial z]\}$ is an i.i.d. sequence with mean zero (which can be shown using integration by parts) and finite variance,

$$
\text{E} \left\{ \frac{1}{n^{1/\alpha_0}} \sum_{t=p+1}^{n} \left[ 1 + Z_t \frac{\partial \ln f(Z_t; \boldsymbol{\tau}_0)}{\partial z} \right] \right\}^2
$$

$$= \frac{1}{n^{2/\alpha_0}} \sum_{t=p+1}^{n} \mathrm{E}\left\{1 + Z_t \frac{\partial \ln f(Z_t; \boldsymbol{\tau}_0)}{\partial z}\right\}^2,$$

which is $o(1)$. Therefore, because $\sup_{\mathbf{u}\in[-T,T]^p}|c_0(\mathbf{u})| < \infty$, (A.22) is $o_p(1)$. Next, (A.23) equals $\sup_{\mathbf{u}\in[-T,T]^p}|(n-p)[n^{-1/\alpha_0} u_p \theta_{0p}^{-1} - \ln|1 + n^{-1/\alpha_0} u_p \theta_{0p}^{-1}|]|$, which is $o(1)$. And finally, (A.24) is bounded above by

$$\sup_{z\in\mathbb{R}}\left|z^2 \frac{\partial^2 \ln f(z; \boldsymbol{\tau}_0)}{\partial z^2}\right| \sup_{\mathbf{u}\in[-T,T]^p} \frac{c_0^2(\mathbf{u})}{2n^{2/\alpha_0}} \sum_{t=p+1}^{n} \left[\frac{Z_t}{Z_{t,n}^\dagger(\mathbf{u})}\right]^2$$

$$\leq \sup_{z\in\mathbb{R}}\left|z^2 \frac{\partial^2 \ln f(z; \boldsymbol{\tau}_0)}{\partial z^2}\right| n^{1-2/\alpha_0} \sup_{\mathbf{u}\in[-T,T]^p} \frac{c_0^2(\mathbf{u})}{2}\left[1 - \frac{|u_p \theta_{0p}^{-1}|}{n^{1/\alpha_0}}\right]^{-2},$$

which, since $\sup_{z\in\mathbb{R}}|z^2[\partial^2 \ln f(z; \boldsymbol{\tau}_0)/\partial z^2]| < \infty$, is also $o(1)$. $\square$

LEMMA A.8. *For any fixed* $\mathbf{u}\in\mathbb{R}^p$ *and* $\mathbf{v}\in\mathbb{R}^4$, $(W_n^\dagger(\mathbf{u}), T_n(\mathbf{v}))' \xrightarrow{\mathcal{L}} (W(\mathbf{u}), \mathbf{v}'\mathbf{N})'$ *on* $\mathbb{R}^2$ *as* $n\to\infty$, *with* $W(\mathbf{u})$ *and* $\mathbf{v}'\mathbf{N}$ *independent.* [$W_n^\dagger(\cdot)$, $T_n(\cdot)$ *and* $W(\cdot)$ *were defined in equations (3.8), (3.9) and (3.4), respectively, and, from Theorem 3.3,* $\mathbf{N}\sim N(\mathbf{0}, \mathbf{I}(\boldsymbol{\tau}_0))$.]

Before proving this result, we introduce some notation and three additional lemmas which will be used in the proof. First, define a set function $\varepsilon_x(\cdot)$ as follows: $\varepsilon_x(A) = I\{x\in A\}$, and, for $m\geq 1$, let

$$\mathbf{e}_1 = (\underbrace{0,\ldots,0}_{m\text{ times}}, 1, \underbrace{0,\ldots,0}_{m-1\text{ times}}),\ldots,\mathbf{e}_m = (\underbrace{0,\ldots,0}_{2m-1\text{ times}}, 1)$$

and

$$\mathbf{e}_{-1} = (\underbrace{0,\ldots,0}_{m-1\text{ times}}, 1, \underbrace{0,\ldots,0}_{m\text{ times}}),\ldots,\mathbf{e}_{-m} = (1, \underbrace{0,\ldots,0}_{2m-1\text{ times}}).$$

Now define

$$S_{m,n}(\cdot) = \sum_{t=p+1}^{n} \varepsilon_{(Z_t, [\tilde{c}(\alpha_0)]^{-1/\alpha_0} \sigma_0^{-1} n^{-1/\alpha_0}(Z_{t+m},\ldots,Z_{t+1},Z_{t-1},\ldots,Z_{t-m}))}(\cdot)$$

and

$$S_m(\cdot) = \sum_{k=1}^{\infty} \sum_{j=1}^{m} (\varepsilon_{(Z_{k,-j}, \mathbf{e}_{-j}\delta_k \Gamma_k^{-1/\alpha_0})}(\cdot) + \varepsilon_{(Z_{k,j}, \mathbf{e}_j\delta_k \Gamma_k^{-1/\alpha_0})}(\cdot)).$$

By the following lemma, $S_{m,n}(\cdot)$ can converge in distribution to $S_m(\cdot)$.

LEMMA A.9. *For any fixed relatively compact subset $A$ of $\mathbb{R}\times(\overline{\mathbb{R}}^{2m}\setminus\{\mathbf{0}\})$ (a subset $A$ for which the closure $\overline{A}$ is compact; note that a compact* *subset of* $\overline{\mathbb{R}}^{2m} \setminus \{\mathbf{0}\} = [-\infty, \infty]^{2m} \setminus \{\mathbf{0}\}$ *is closed and bounded away from the origin) of the form*

$$
\begin{aligned}
(A.25) \quad A = (a_0, b_0] &\times (a_{-m}, b_{-m}] \times \cdots \times (a_{-1}, b_{-1}] \\
&\times (a_1, b_1] \times \cdots \times (a_m, b_m], \qquad a_j, b_j \neq 0 \ \forall |j| \in \{1, \ldots, m\}
\end{aligned}
$$

*and for any fixed* $\mathbf{v} \in \mathbb{R}^4$, $(S_{m,n}(A), T_n(\mathbf{v}))' \xrightarrow{\mathcal{L}} (S_m(A), \mathbf{v}'\mathbf{N})'$ *on* $\mathbb{R}^2$ *as* $n \to \infty$, *with* $S_m(A)$ *and* $\mathbf{v}'\mathbf{N}$ *independent.*

PROOF. Let $\lambda_1, \lambda_2 \in \mathbb{R}$. Following Theorem 3 on page 37 of Rosenblatt [34], this Lemma holds if $\text{cum}_k(\lambda_1 S_{m,n}(A) + \lambda_2 T_n(\mathbf{v})) \to \text{cum}_k(\lambda_1 S_m(A) + \lambda_2 \mathbf{v}'\mathbf{N})$ for all $k \geq 1$, where $\text{cum}_k(X)$ is the $k$th-order joint cumulant of the random variable $X$. So,

$$
\text{cum}_k(X) = \text{cum}(\underbrace{X, \ldots, X}_{k \text{ times}}).
$$

Note that since $S_m(A)$ and $\mathbf{v}'\mathbf{N}$ are independent, $\text{cum}_k(\lambda_1 S_m(A) + \lambda_2 \mathbf{v}'\mathbf{N}) = \lambda_1^k \text{cum}_k(S_m(A)) + \lambda_2^k \text{cum}_k(\mathbf{v}'\mathbf{N})$.

Fix $k \geq 1$ and denote the $k$th-order joint cumulant of $i$ $X$s and $j$ $Y$s $(i + j = k)$ as $\text{cum}_{i,j}(X, Y)$. So,

$$
\text{cum}_{i,j}(X, Y) = \text{cum}(\underbrace{X, \ldots, X}_{i \text{ times}}, \underbrace{Y, \ldots, Y}_{j \text{ times}}).
$$

Then, by linearity,

$$
\begin{aligned}
\text{cum}_k(\lambda_1 S_{m,n}(A) + \lambda_2 T_n(\mathbf{v})) = {}& \lambda_1^k \text{cum}_k(S_{m,n}(A)) + \lambda_2^k \text{cum}_k(T_n(\mathbf{v})) \\
&+ \sum_{j=1}^{k-1} \binom{k}{j} \lambda_1^j \lambda_2^{k-j} \text{cum}_{j,k-j}(S_{m,n}(A), T_n(\mathbf{v})).
\end{aligned}
$$

Also by linearity, for $j \in \{1, \ldots, k-1\}$,

$$
\begin{aligned}
(A.26) \quad &\text{cum}_{j,k-j}(S_{m,n}(A), T_n(\mathbf{v})) \\
&= \sum_{t_1=p+1}^{n} \cdots \sum_{t_k=p+1}^{n} \text{cum}(V_{t_1,n}, \ldots, V_{t_j,n}, W_{t_{j+1},n}, \ldots, W_{t_k,n}),
\end{aligned}
$$

where $V_{t,n} := \varepsilon_{(Z_t, [\tilde{c}(\alpha_0)]^{-1/\alpha_0} \sigma_0^{-1} n^{-1/\alpha_0}(Z_{t+m}, \ldots, Z_{t+1}, Z_{t-1}, \ldots, Z_{t-m}))}(A)$ and $W_{t,n} := n^{-1/2} \mathbf{v}' \partial \ln f(Z_t; \boldsymbol{\tau}_0)/\partial \boldsymbol{\tau}$. Due to the limited dependence between the variables $\{V_{t,n}\}_{t=p+1}^{n}$, $\{W_{t,n}\}_{t=p+1}^{n}$, equation (A.26) equals

$$
\begin{aligned}
(A.27) \quad &\sum_{t_1=p+1}^{n} \sum_{|t_2-t_1| \leq 2jm} \cdots \sum_{|t_j-t_1| \leq 2jm} \sum_{|t_{j+1}-t_1| \leq (2j+1)m} \\
&\cdots \sum_{|t_k-t_1| \leq (2j+1)m} \text{cum}(V_{t_1,n}, \ldots, V_{t_j,n}, W_{t_{j+1},n}, \ldots, W_{t_k,n});
\end{aligned}
$$

this sum is made up of $(n-p)(4jm+1)^{j-1}([4j+2]m+1)^{k-j}$ terms. Therefore, since $|V_{t,n}| \leq 1$, and $\mathrm{E}|W_{t,n}|^{\ell} < \infty$ for all $\ell \geq 1$ and all $n$, (A.27) is $o(1)$ if $k-j \geq 3$ [as a result of the scaling by $n^{-(k-j)/2}$]. Equation (A.27) is also $o(1)$ for $k-j \in \{1,2\}$ if $n\mathrm{E}|V_{t_1,n}W_{t_2,n}| = o(1)$ and $n\mathrm{E}|V_{t_1,n}W_{t_2,n}W_{t_3,n}| = o(1)$ for any $t_1,t_2,t_3$. We will show the limit is zero in one case; convergence to zero can be established similarly in all other cases.

Since $A$ is a relatively compact subset of $\mathbb{R} \times (\overline{\mathbb{R}}^{2m} \setminus \{\mathbf{0}\})$, at least one of the intervals $(a_{-m}, b_{-m}], \ldots, (a_{-1}, b_{-1}], (a_1, b_1], \ldots, (a_m, b_m]$ does not contain zero. We assume $(a_{-1}, b_{-1}]$ does not contain zero and show that $n\mathrm{E}|V_{1,n}W_{2,n}| = o(1)$. First, from (2.5)–(2.10), there exist constants $C_v, D_v < \infty$ such that $|\mathbf{v}'\partial \ln f(z; \boldsymbol{\tau}_0)/\partial \boldsymbol{\tau}| \leq C_v + D_v|z|^{\alpha_0/4}$ $\quad \forall z \in \mathbb{R}$. Hence, because $V_{1,n} = \varepsilon_{(Z_1, [\tilde{c}(\alpha_0)]^{-1/\alpha_0}\sigma_0^{-1}n^{-1/\alpha_0}(Z_{1+m}, \ldots, Z_2, Z_0, \ldots, Z_{1-m}))}(A)$ and $W_{2,n} = n^{-1/2}\mathbf{v}'\partial \ln f(Z_2; \boldsymbol{\tau}_0)/\partial \boldsymbol{\tau}$,

$$n\mathrm{E}|V_{1,n}W_{2,n}|$$
$$\leq n^{1/2}\mathrm{E}\bigg|I\{[\tilde{c}(\alpha_0)]^{-1/\alpha_0}\sigma_0^{-1}n^{-1/\alpha_0}Z_2 \in (a_{-1}, b_{-1}]\}\bigg(\mathbf{v}'\frac{\partial \ln f(Z_2; \boldsymbol{\tau}_0)}{\partial \boldsymbol{\tau}}\bigg)\bigg|$$
$$\leq C_v n^{1/2}\mathrm{P}(|Z_2| \geq n^{1/\alpha_0}\zeta) + D_v n^{1/2}\mathrm{E}\{|Z_2|^{\alpha_0/4}I\{|Z_2| \geq n^{1/\alpha_0}\zeta\}\},$$

where $\zeta := [\tilde{c}(\alpha_0)]^{1/\alpha_0}\sigma_0 \min\{|a_{-1}|, |b_{-1}|\}$. By (2.3), since $\zeta > 0$, $n^{1/2}\mathrm{P}(|Z_2| \geq n^{1/\alpha_0}\zeta) \to 0$, and, using Karamata's theorem (see, e.g., Feller [19], page 283), $n^{1/2}\mathrm{E}\{|Z_2|^{\alpha_0/4}I\{|Z_2| \geq n^{1/\alpha_0}\zeta\}\} \leq (constant)n^{1/2}(n^{1/\alpha_0}\zeta)^{\alpha_0/4}\mathrm{P}(|Z_2| \geq n^{1/\alpha_0}\zeta)$, which is $o(1)$ by (2.3).

It has therefore been established that $\mathrm{cum}_k(\lambda_1 S_{m,n}(A) + \lambda_2 T_n(\mathbf{v})) = \lambda_1^k \times \mathrm{cum}_k(S_{m,n}(A)) + \lambda_2^k \mathrm{cum}_k(T_n(\mathbf{v})) + o(1)$ for arbitrary $k \geq 1$. Following the Proof of Lemma 16 in Calder [7], it can be shown that $\mathrm{cum}_k(S_{m,n}(A)) \to \mathrm{cum}_k(S_m(A))$. Note that, from Davis and Resnick [13], $S_{m,n}(A) \xrightarrow{\mathcal{L}} S_m(A)$ on $\mathbb{R}$ and $S_m(A)$ is a Poisson random variable, so all cumulants are finite. It is relatively straightforward to show that $\mathrm{cum}_k(T_n(\mathbf{v})) \to \mathrm{cum}_k(\mathbf{v}'\mathbf{N})$ (see the Proof of Lemma 16 in [7] for details), which is not surprising since $T_n(\mathbf{v}) \xrightarrow{\mathcal{L}} \mathbf{v}'\mathbf{N}$ on $\mathbb{R}$ by the central limit theorem. Consequently, $\mathrm{cum}_k(\lambda_1 S_{m,n}(A) + \lambda_2 T_n(\mathbf{v})) \to \lambda_1^k \mathrm{cum}_k(S_m(A)) + \lambda_2^k \mathrm{cum}_k(\mathbf{v}'\mathbf{N})$, and the proof is complete. $\quad \square$

LEMMA A.10.   Let $U_{t,n}^-(\mathbf{u}) = n^{-1/\alpha_0}\sum_{j=1}^{\infty} c_{-j}(\mathbf{u})Z_{t+j}$, $U_{t,n}^+(\mathbf{u}) = n^{-1/\alpha_0} \times \sum_{j=1}^{\infty} c_j(\mathbf{u})Z_{t-j}$ and $I_{t,n}^{\lambda,\lambda,M} = I\{|Z_t| \leq M\}I\{(|U_{t,n}^-(\mathbf{u})| > \lambda) \cup (|U_{t,n}^+(\mathbf{u})| > \lambda)\}$. For any fixed $\mathbf{u} \in \mathbb{R}^p$ and any $\kappa > 0$, $\lim_{\lambda \to 0^+} \lim_{M \to \infty} \limsup_{n \to \infty}$ of

$$\mathrm{P}\bigg(\bigg|\sum_{t=p+1}^{n}\{[\ln f(Z_t + U_{t,n}^-(\mathbf{u}) + U_{t,n}^+(\mathbf{u}); \boldsymbol{\tau}_0)$$
$$(A.28)$$
$$- \ln f(Z_t; \boldsymbol{\tau}_0)][1 - I_{t,n}^{\lambda,\lambda,M}]\}\bigg| > \kappa\bigg)$$

*is zero.*

PROOF.   Note that, for any $t \in \{p+1, \dots, n\}$ and any $n$, $\ln f(Z_t + U_{t,n}^-(\mathbf{u}) + U_{t,n}^+(\mathbf{u}); \boldsymbol{\tau}_0) - \ln f(Z_t; \boldsymbol{\tau}_0) = [U_{t,n}^-(\mathbf{u}) + U_{t,n}^+(\mathbf{u})][\partial \ln f(Z_t; \boldsymbol{\tau}_0)/\partial z] + [U_{t,n}^-(\mathbf{u}) + U_{t,n}^+(\mathbf{u})]^2 [\partial^2 \ln f(Z_{t,n}^*; \boldsymbol{\tau}_0)/\partial z^2]/2$, where $Z_{t,n}^*$ lies between $Z_t$ and $Z_t + U_{t,n}^-(\mathbf{u}) + U_{t,n}^+(\mathbf{u})$. Note also that

$$1 - I_{t,n}^{\lambda,\lambda,M} = I\{|U_{t,n}^-(\mathbf{u})| \le \lambda\} I\{|U_{t,n}^+(\mathbf{u})| \le \lambda\} + I\{|Z_t| > M\} I\{|U_{t,n}^-(\mathbf{u})| > \lambda\}$$
$$+ I\{|Z_t| > M\} I\{|U_{t,n}^+(\mathbf{u})| > \lambda\}$$
$$- I\{|Z_t| > M\} I\{|U_{t,n}^-(\mathbf{u})| > \lambda\} I\{|U_{t,n}^+(\mathbf{u})| > \lambda\}.$$

Consequently, (A.28) is bounded above by

$$P\left(\left|\sum_{t=p+1}^n \left[ (U_{t,n}^-(\mathbf{u}) + U_{t,n}^+(\mathbf{u})) \frac{\partial \ln f(Z_t; \boldsymbol{\tau}_0)}{\partial z} I\{|U_{t,n}^-(\mathbf{u})| \le \lambda\} \right. \right. \right.$$
$$\left. \left. \left. \times I\{|U_{t,n}^+(\mathbf{u})| \le \lambda\} \right] \right| > \frac{\kappa}{5} \right)$$

$$+ P\left( \sup_{z \in \mathbb{R}} \left| \frac{\partial^2 \ln f(z; \boldsymbol{\tau}_0)}{\partial z^2} \right| \right.$$
$$\left. \times \sum_{t=p+1}^n [U_{t,n}^-(\mathbf{u}) + U_{t,n}^+(\mathbf{u})]^2 I\{|U_{t,n}^-(\mathbf{u})| \le \lambda\} I\{|U_{t,n}^+(\mathbf{u})| \le \lambda\}] > \frac{\kappa}{5} \right)$$

$$+ P\left( \bigcup_{t=p+1}^n \{(|Z_t| > M) \cap (|U_{t,n}^-(\mathbf{u})| > \lambda)\} \right)$$

$$+ 2P\left( \bigcup_{t=p+1}^n \{(|Z_t| > M) \cap (|U_{t,n}^+(\mathbf{u})| > \lambda)\} \right).$$

The proof of Proposition A.2(a)–(c) in Davis, Knight and Liu [12] can be used to show that $\lim_{\lambda \to 0^+} \lim_{M \to \infty} \lim \sup_{n \to \infty}$ of each of the four summands is zero.   □

LEMMA A.11.   *Let* $I_{k,j}^{\lambda,M} = I\{|Z_{k,j}| \le M\} I\{|[\tilde{c}(\alpha_0)]^{1/\alpha_0} \sigma_0 c_j(\mathbf{u}) \delta_k \Gamma_k^{-1/\alpha_0}| > \lambda\}$. *For any fixed* $\mathbf{u} \in \mathbb{R}^p$,

$$\sum_{k=1}^\infty \sum_{j \ne 0} [\{\ln f(Z_{k,j} + [\tilde{c}(\alpha_0)]^{1/\alpha_0} \sigma_0 c_j(\mathbf{u}) \delta_k \Gamma_k^{-1/\alpha_0}; \boldsymbol{\tau}_0)$$
(A.29)
$$- \ln f(Z_{k,j}; \boldsymbol{\tau}_0)\}(1 - I_{k,j}^{\lambda,M})]$$

*converges in probability to zero as* $\lambda \to 0^+$ *and* $M \to \infty$.

PROOF. The absolute value of (A.29) is bounded above by $[\tilde{c}(\alpha_0)]^{1/\alpha_0} \times \sigma_0 \sup_{z \in \mathbb{R}} |\partial \ln f(z; \boldsymbol{\tau}_0)/\partial z| \sum_{k=1}^{\infty} \Gamma_k^{-1/\alpha_0} \sum_{j \neq 0} |c_j(\mathbf{u})|$. If $\alpha_0 < 1$, $\sum_{k=1}^{\infty} \Gamma_k^{-1/\alpha_0} < \infty$ a.s., since $\mathrm{E}\{\Gamma_k^{-1/\alpha_0}\} = O(k^{-1/\alpha_0})$ for $k > 1/\alpha_0$. Thus, the result holds if $\alpha_0 < 1$.

For $\alpha_0 \geq 1$, the proof of this lemma is similar to the Proof of Lemma A.10. We omit the details. $\square$

We now use Lemmas A.9–A.11 to prove Lemma A.8.

PROOF OF LEMMA A.8. By Lemma A.9, for any relatively compact subset $A$ of $\mathbb{R} \times (\overline{\mathbb{R}}^{2m} \setminus \{\mathbf{0}\})$ of the form (A.25) and any $\mathbf{v} \in \mathbb{R}^4$, $(S_{m,n}(A), T_n(\mathbf{v}))' \xrightarrow{\mathcal{L}} (S_m(A), \mathbf{v}'\mathbf{N})'$ on $\mathbb{R}^2$, with $S_m(A)$ and $\mathbf{v}'\mathbf{N}$ independent. It can be shown similarly that, for any $\ell \geq 1$ and any relatively compact subsets $A_1, \ldots, A_\ell$ of $\mathbb{R} \times (\overline{\mathbb{R}}^{2m} \setminus \{\mathbf{0}\})$ of the form (A.25),

$$(A.30) \quad (S_{m,n}(A_1), \ldots, S_{m,n}(A_\ell), T_n(\mathbf{v}))' \xrightarrow{\mathcal{L}} (S_m(A_1), \ldots, S_m(A_\ell), \mathbf{v}'\mathbf{N})'$$

on $\mathbb{R}^{\ell+1}$, with $(S_m(A_1), \ldots, S_m(A_\ell))'$ and $\mathbf{v}'\mathbf{N}$ independent. Now, for fixed $\mathbf{u} \in \mathbb{R}^p$, let $\tilde{S}_n(\cdot) = \sum_{t=p+1}^n \varepsilon_{(Z_t, U_{t,n}^-(\mathbf{u}), U_{t,n}^+(\mathbf{u}))}(\cdot)$, with $U_{t,n}^-(\mathbf{u}) = n^{-1/\alpha_0} \times \sum_{j=1}^{\infty} c_{-j}(\mathbf{u}) Z_{t+j}$ and $U_{t,n}^+(\mathbf{u}) = n^{-1/\alpha_0} \sum_{j=1}^{\infty} c_j(\mathbf{u}) Z_{t-j}$, and let

$$\tilde{S}(\cdot) = \sum_{k=1}^{\infty} \sum_{j=1}^{\infty} (\varepsilon_{(Z_{k,-j}, [\tilde{c}(\alpha_0)]^{1/\alpha_0} \sigma_0 c_{-j}(\mathbf{u}) \delta_k \Gamma_k^{-1/\alpha_0}, 0)}(\cdot)$$
$$+ \varepsilon_{(Z_{k,j}, 0, [\tilde{c}(\alpha_0)]^{1/\alpha_0} \sigma_0 c_j(\mathbf{u}) \delta_k \Gamma_k^{-1/\alpha_0})}(\cdot)).$$

Following the proof of Theorem 2.4 in Davis and Resnick [13], using (A.30), the mapping

$$(z_t, z_{t+m}, \ldots, z_{t+1}, z_{t-1}, \ldots, z_{t-m})$$
$$\to \left(z_t, [\tilde{c}(\alpha_0)]^{1/\alpha_0} \sigma_0 \sum_{j=1}^{m} c_{-j}(\mathbf{u}) z_{t+j}, [\tilde{c}(\alpha_0)]^{1/\alpha_0} \sigma_0 \sum_{j=1}^{m} c_j(\mathbf{u}) z_{t-j}\right),$$

and by letting $m \to \infty$, it can be shown that

$$(A.31) \quad (\tilde{S}_n(\tilde{A}_1), \ldots, \tilde{S}_n(\tilde{A}_\ell), T_n(\mathbf{v}))' \xrightarrow{\mathcal{L}} (\tilde{S}(\tilde{A}_1), \ldots, \tilde{S}(\tilde{A}_\ell), \mathbf{v}'\mathbf{N})'$$

on $\mathbb{R}^{\ell+1}$, with $(\tilde{S}(\tilde{A}_1), \ldots, \tilde{S}(\tilde{A}_\ell))'$ and $\mathbf{v}'\mathbf{N}$ independent, for any relatively compact subsets $\tilde{A}_1, \ldots, \tilde{A}_\ell$ of $\mathbb{R} \times (\overline{\mathbb{R}}^2 \setminus \{\mathbf{0}\})$.

Since $(\tilde{S}_n(\tilde{A}_1), \ldots, \tilde{S}_n(\tilde{A}_\ell))' \xrightarrow{\mathcal{L}} (\tilde{S}(\tilde{A}_1), \ldots, \tilde{S}(\tilde{A}_\ell))'$ on $\mathbb{R}^\ell$ for arbitrary $\ell \geq 1$ and arbitrary, relatively compact subsets $\tilde{A}_1, \ldots, \tilde{A}_\ell$ of $\mathbb{R} \times (\overline{\mathbb{R}}^2 \setminus \{\mathbf{0}\})$,

$$\sum_{t=p+1}^{n} \tilde{g}(Z_t, U_{t,n}^-(\mathbf{u}), U_{t,n}^+(\mathbf{u}))$$

$$(A.32) \qquad \xrightarrow{\mathcal{L}} \sum_{k=1}^{\infty} \sum_{j=1}^{\infty} (\tilde{g}(Z_{k,-j}, [\tilde{c}(\alpha_0)]^{1/\alpha_0} \sigma_0 c_{-j}(\mathbf{u}) \delta_k \Gamma_k^{-1/\alpha_0}, 0)$$

$$+ \tilde{g}(Z_{k,j}, 0, [\tilde{c}(\alpha_0)]^{1/\alpha_0} \sigma_0 c_j(\mathbf{u}) \delta_k \Gamma_k^{-1/\alpha_0}))$$

on $\mathbb{R}$ for any continuous function $\tilde{g}$ on $\mathbb{R} \times (\overline{\mathbb{R}}^2 \setminus \{\mathbf{0}\})$ with compact support (see Davis and Resnick [13]). Because it is almost everywhere continuous on $\mathbb{R} \times (\overline{\mathbb{R}}^2 \setminus \{\mathbf{0}\})$ with compact support, we will use $\tilde{g}(x, y, z) = [\ln f(x + y + z; \boldsymbol{\tau}_0) - \ln f(x; \boldsymbol{\tau}_0)] I\{|x| \le M\} I\{(|y| > \lambda) \cup (|z| > \lambda)\}$, where $M, \lambda > 0$. By Lemma A.10, for any $\kappa > 0$, $\lim_{\lambda \to 0^+} \lim_{M \to \infty} \limsup_{n \to \infty} \mathrm{P}(|W_n^\dagger(\mathbf{u}) - \sum_{t=p+1}^{n} \tilde{g}(Z_t, U_{t,n}^-(\mathbf{u}), U_{t,n}^+(\mathbf{u}))| > \kappa) = 0$ and, by Lemma A.11,

$$\sum_{k=1}^{\infty} \sum_{j=1}^{\infty} (\tilde{g}(Z_{k,-j}, [\tilde{c}(\alpha_0)]^{1/\alpha_0} \sigma_0 c_{-j}(\mathbf{u}) \delta_k \Gamma_k^{-1/\alpha_0}, 0)$$

$$+ \tilde{g}(Z_{k,j}, 0, [\tilde{c}(\alpha_0)]^{1/\alpha_0} \sigma_0 c_j(\mathbf{u}) \delta_k \Gamma_k^{-1/\alpha_0})) \xrightarrow{P} W(\mathbf{u})$$

as $\lambda \to 0^+$ and $M \to \infty$ [$W_n^\dagger(\cdot)$ and $W(\cdot)$ were defined in equations (3.8) and (3.4), resp.]. Therefore, by Theorem 3.2 in Billingsley [2], it follows from (A.32) that $W_n^\dagger(\mathbf{u}) \xrightarrow{\mathcal{L}} W(\mathbf{u})$ on $\mathbb{R}$ for fixed $\mathbf{u} \in \mathbb{R}^p$, and consequently the result of this lemma follows from (A.31). $\square$

LEMMA A.12. *For any $T > 0$ and any $\kappa > 0$,*

$$(A.33) \quad \lim_{\epsilon \to 0^+} \limsup_{n \to \infty} P\left( \sup_{\|\mathbf{u}\|, \|\mathbf{v}\| \le T, \|\mathbf{u} - \mathbf{v}\| \le \epsilon} |W_n^\dagger(\mathbf{u}) - W_n^\dagger(\mathbf{v})| > \kappa \right) = 0.$$

*[$W_n^\dagger(\cdot)$ was defined in equation (3.8).]*

PROOF. For $\mathbf{u}, \mathbf{v} \in \mathbb{R}^p$,

$$|W_n^\dagger(\mathbf{u}) - W_n^\dagger(\mathbf{v})| = \left| \frac{1}{n^{1/\alpha_0}} \sum_{t=p+1}^{n} \left( \sum_{j \ne 0} c_j(\mathbf{u} - \mathbf{v}) Z_{t-j} \right) \frac{\partial \ln f(Z_{t,n}^*(\mathbf{u}, \mathbf{v}); \boldsymbol{\tau}_0)}{\partial z} \right|$$

$$\le \left| \frac{1}{n^{1/\alpha_0}} \sum_{t=p+1}^{n} \left( \sum_{j \ne 0} c_j(\mathbf{u} - \mathbf{v}) Z_{t-j} \right) \frac{\partial \ln f(Z_t; \boldsymbol{\tau}_0)}{\partial z} \right|$$

$$+ \left( \sup_{z \in \mathbb{R}} \left| \frac{\partial^2 \ln f(z; \boldsymbol{\tau}_0)}{\partial z^2} \right| \right)$$

$$\times \frac{1}{n^{2/\alpha_0}} \sum_{t=p+1}^{n} \left| \sum_{j \ne 0} c_j(\mathbf{u} - \mathbf{v}) Z_{t-j} \right|$$

$$\times \left[ \left| \sum_{j \ne 0} c_j(\mathbf{u}) Z_{t-j} \right| + \left| \sum_{j \ne 0} c_j(\mathbf{v}) Z_{t-j} \right| \right],$$

where $Z_{t,n}^*(\mathbf{u}, \mathbf{v})$ lies between $Z_t + n^{-1/\alpha_0} \sum_{j\neq 0} c_j(\mathbf{u}) Z_{t-j}$ and $Z_t + n^{-1/\alpha_0} \times \sum_{j\neq 0} c_j(\mathbf{v}) Z_{t-j}$. Following the Proof of Theorem 2.1 in Davis, Knight and Liu [12] (see page 154), if $\{\tilde{\pi}_j\}_{j\neq 0}$ is a geometrically decaying sequence, then it can be shown that $n^{-1/\alpha_0} \sum_{t=p+1}^n (\sum_{j\neq 0} \tilde{\pi}_j Z_{t-j})[\partial \ln f(Z_t; \boldsymbol{\tau}_0)/\partial z] = O_p(1)$ and $n^{-2/\alpha_0} \sum_{t=p+1}^n (\sum_{j\neq 0} |\tilde{\pi}_j Z_{t-j}|)^2 = O_p(1)$. Therefore, by (A.12) and because $c_j(\mathbf{u})$ is linear in $\mathbf{u}$ for all $j$, (A.33) holds. $\quad\square$

LEMMA A.13. *If, as $n \to \infty$, $m_n \to \infty$ with $m_n/n \to 0$, then for any $T > 0$ and any $\kappa > 0$,*

$$(A.34) \qquad P\Big(\sup_{\|\mathbf{u}\| \leq T} |\tilde{W}_{m_n}(\mathbf{u}) - \tilde{W}_{m_n}^\dagger(\mathbf{u})| > \kappa | X_1, \ldots, X_n\Big) \xrightarrow{P} 0.$$

*[$\tilde{W}_{m_n}^\dagger(\cdot)$ and $\tilde{W}_{m_n}(\cdot)$ were defined in equations (3.14) and (3.15).]*

PROOF. Choose arbitrary $T, \kappa > 0$, and let the sequence $\{\hat{\psi}_j\}_{j=-\infty}^\infty$ contain the coefficients in the Laurent series expansion of $1/[\hat{\theta}_{\mathrm{ML}}^\dagger(z)\hat{\theta}_{\mathrm{ML}}^*(z)]$. From (3.11), for $t \in \{1, \ldots, m_n\}$, $\hat{\theta}_{\mathrm{ML}}^\dagger(B)\hat{\theta}_{\mathrm{ML}}^*(B)X_t^* = Z_t^*$, and so $X_t^* = \sum_{j=-\infty}^\infty \hat{\psi}_j Z_{t-j}^*$. From Brockwell and Davis [6] (see Chapter 3), there exist $C_2 > 0$, $0 < D_2 < 1$ and a sufficiently small $\delta > 0$ such that, whenever $\|\hat{\boldsymbol{\theta}}_{\mathrm{ML}} - \boldsymbol{\theta}_0\| < \delta$, $|\hat{\psi}_j| \leq C_2 D_2^{|j|}$ and also $\sup_{\|\mathbf{u}\| \leq T} |\hat{c}_j(\mathbf{u})| \leq C_2 D_2^{|j|}$ for all $j \in \{\ldots, -1, 0, 1, \ldots\}$ [the $\hat{c}_j(\mathbf{u})$s were defined in (3.13)]. Now observe that the left-hand side of (A.34) is bounded above by

$$
\begin{aligned}
(A.35) \quad & P\Big(\sup_{\|\mathbf{u}\| \leq T} |\tilde{W}_{m_n}(\mathbf{u}) - \tilde{W}_{m_n}^\dagger(\mathbf{u})| > \kappa | X_1, \ldots, X_n\Big) I\{\|\hat{\boldsymbol{\theta}}_{\mathrm{ML}} - \boldsymbol{\theta}_0\| < \delta\} \\
& + I\{\|\hat{\boldsymbol{\theta}}_{\mathrm{ML}} - \boldsymbol{\theta}_0\| \geq \delta\},
\end{aligned}
$$

and that $I\{\|\hat{\boldsymbol{\theta}}_{\mathrm{ML}} - \boldsymbol{\theta}_0\| \geq \delta\}$ is $o_p(1)$ since $\hat{\boldsymbol{\theta}}_{\mathrm{ML}} \xrightarrow{P} \boldsymbol{\theta}_0$. For $\mathbf{u} = (u_1, \ldots, u_p)' \in \mathbb{R}^p$,

$$
\begin{aligned}
& \tilde{W}_{m_n}(\mathbf{u}) - \tilde{W}_{m_n}^\dagger(\mathbf{u}) \\
& \quad = \sum_{t=p+1}^{m_n} \Big[\ln f\Big(Z_t^*\Big(\hat{\boldsymbol{\theta}}_{\mathrm{ML}} + \frac{\mathbf{u}}{m_n^{1/\alpha_0}}, s_0\Big); \hat{\boldsymbol{\tau}}_{\mathrm{ML}}\Big) - \ln f(Z_t^*; \hat{\boldsymbol{\tau}}_{\mathrm{ML}})\Big] \\
& \qquad - \sum_{t=p+1}^{m_n} \Big[\ln f\Big(Z_t^* + m_n^{-1/\alpha_0} \sum_{j\neq 0} \hat{c}_j(\mathbf{u}) Z_{t-j}^*; \boldsymbol{\tau}_0\Big) - \ln f(Z_t^*; \boldsymbol{\tau}_0)\Big] \\
& \qquad + (m_n - p) \ln\Big|\frac{\hat{\theta}_{p,\mathrm{ML}} + m_n^{-1/\alpha_0} u_p}{\hat{\theta}_{p,\mathrm{ML}}}\Big| I\{s_0 > 0\},
\end{aligned}
$$

and so, using arguments similar to those given in the proofs of Lemmas A.4–A.7, it can be shown that the first summand of (A.35) is also $o_p(1)$ if, for any $\epsilon > 0$,

$$
\text{(A.36)} \qquad \mathrm{P}\left( \frac{C_2}{m_n^{2/\alpha_0}} \sum_{t=p+1}^{m_n} \sum_{j=-\infty}^{\infty} D_2^{|j|} |Z_{t-j}^*| > \epsilon \,|\, X_1, \ldots, X_n \right),
$$

$$
\text{(A.37)} \qquad
\begin{aligned}
&\mathrm{P}\Bigg( \left[ \sup_{i \in \{1,\ldots,4\}} |\hat{\tau}_{i,\mathrm{ML}} - \tau_{0i}| \right] \frac{C_2}{m_n^{1/\alpha_0}} \\
&\qquad \times \sum_{t=p+1}^{m_n} \sum_{j=-\infty}^{\infty} D_2^{|j|} |Z_{t-j}^*| > \epsilon \,|\, X_1, \ldots, X_n \Bigg),
\end{aligned}
$$

$$
\text{(A.38)} \qquad \mathrm{P}\left( \frac{C_2}{m_n^{2/\alpha_0}} \sum_{t=p+1}^{m_n} |Z_t^*| \sum_{j \neq 0} D_2^{|j|} |Z_{t-j}^*| > \epsilon \,|\, X_1, \ldots, X_n \right),
$$

$$
\text{(A.39)} \qquad \mathrm{P}\left( \frac{C_2}{m_n^{3/\alpha_0}} \sum_{t=p+1}^{m_n} (Z_t^*)^2 \sum_{j \neq 0} D_2^{|j|} |Z_{t-j}^*| > \epsilon \,|\, X_1, \ldots, X_n \right)
$$

and

$$
\text{(A.40)} \quad \mathrm{P}\left( \left\{ \frac{1}{m_n^{1/\alpha_0}} \sum_{t=p+1}^{m_n} \left[ 1 + Z_t^* \frac{\partial \ln f(Z_t^*; \boldsymbol{\tau}_0)}{\partial z} \right] \right\}^2 > \epsilon \,|\, X_1, \ldots, X_n \right)
$$

are all $o_p(1)$. To complete the proof, we show that (A.36) and (A.40) are both $o_p(1)$. Since $n^{1/2}(\hat{\boldsymbol{\tau}}_{\mathrm{ML}} - \boldsymbol{\tau}_0) = O_p(1)$ and $m_n/n \to 0$, using the Proof of Lemma A.5, it can be shown similarly that (A.37) is $o_p(1)$. The Proof of Lemma A.6 can be used to show that (A.38) and (A.39) are $o_p(1)$.

Recall, from the Proof of Lemma A.4, that $\kappa_1 = (3/4)\alpha_0 I\{\alpha_0 \leq 1\} + I\{\alpha_0 > 1\}$. By the Markov inequality, equation (A.36) is bounded above by

$$
\left( \frac{C_2}{\epsilon} \right)^{\kappa_1} m_n^{1 - 2\kappa_1/\alpha_0} \left[ \sum_{j=-\infty}^{\infty} (D_2^{\kappa_1})^{|j|} \right] \mathrm{E}\{ |Z_t^*|^{\kappa_1} \,|\, X_1, \ldots, X_n \};
$$

this is $o_p(1)$ since $m_n^{1 - 2\kappa_1/\alpha_0} \to 0$ and, using $\hat{\boldsymbol{\theta}}_{\mathrm{ML}} \xrightarrow{P} \boldsymbol{\theta}_0$ and $\mathrm{E}|Z_1|^{\kappa_1} < \infty$, it can be shown that $\mathrm{E}\{ |Z_t^*|^{\kappa_1} \,|\, X_1, \ldots, X_n \} = (n-p)^{-1} \sum_{t=p+1}^{n} |Z_t(\hat{\boldsymbol{\theta}}_{\mathrm{ML}}, s_0)|^{\kappa_1}$ is $O_p(1)$.

We now consider (A.40), which is bounded above by

$$
\text{(A.41)} \quad \epsilon^{-1} m_n^{1 - 2/\alpha_0} \mathrm{E}\left\{ \left( 1 + Z_t^* \frac{\partial \ln f(Z_t^*; \boldsymbol{\tau}_0)}{\partial z} \right)^2 \,|\, X_1, \ldots, X_n \right\}
$$

$$
\text{(A.42)} \quad + \epsilon^{-1} \frac{m_n^2 - m_n}{m_n^{2/\alpha_0}} \left[ \mathrm{E}\left\{ 1 + Z_t^* \frac{\partial \ln f(Z_t^*; \boldsymbol{\tau}_0)}{\partial z} \,|\, X_1, \ldots, X_n \right\} \right]^2.
$$

Since $m_n^{1-2/\alpha_0} \to 0$ and, by (2.6), $\sup_{z \in \mathbb{R}} |z[\partial \ln f(z; \boldsymbol{\tau}_0)/\partial z]| < \infty$, (A.41) is $o_p(1)$. Now consider

$$(\text{A.43}) \quad \mathrm{E}\left\{1 + Z_t^* \frac{\partial \ln f(Z_t^*; \boldsymbol{\tau}_0)}{\partial z} | X_1, \ldots, X_n\right\}$$

$$(\text{A.44}) \qquad = \frac{1}{n-p} \sum_{t=p+1}^{n} \left(1 + Z_t \frac{\partial \ln f(Z_t; \boldsymbol{\tau}_0)}{\partial z}\right)$$

$$+ \frac{1}{n-p}\left[\sum_{t=p+1}^{n} Z_t(\hat{\boldsymbol{\theta}}_{\mathrm{ML}}, s_0) \frac{\partial \ln f(Z_t(\hat{\boldsymbol{\theta}}_{\mathrm{ML}}, s_0); \boldsymbol{\tau}_0)}{\partial z}\right.$$

$$(\text{A.45})$$

$$\left. - \sum_{t=p+1}^{n} Z_t \frac{\partial \ln f(Z_t; \boldsymbol{\tau}_0)}{\partial z}\right].$$

By the central limit theorem, (A.44) is $O_p(n^{-1/2})$. In addition, since $Z_t = Z_t(\boldsymbol{\theta}_0, s_0)$, (A.45) equals

$$(\text{A.46}) \quad \frac{(\hat{\boldsymbol{\theta}}_{\mathrm{ML}} - \boldsymbol{\theta}_0)'}{n-p} \sum_{t=p+1}^{n} \left[\frac{\partial \ln f(Z_t(\boldsymbol{\theta}_n^*, s_0); \boldsymbol{\tau}_0)}{\partial z}\right.$$

$$\left. + Z_t(\boldsymbol{\theta}_n^*, s_0) \frac{\partial^2 \ln f(Z_t(\boldsymbol{\theta}_n^*, s_0); \boldsymbol{\tau}_0)}{\partial z^2}\right] \frac{\partial Z_t(\boldsymbol{\theta}_n^*, s_0)}{\partial \boldsymbol{\theta}},$$

with $\boldsymbol{\theta}_n^*$ between $\hat{\boldsymbol{\theta}}_{\mathrm{ML}}$ and $\boldsymbol{\theta}_0$, and, because $\sup_{z \in \mathbb{R}} |[\partial \ln f(z; \boldsymbol{\tau}_0)/\partial z] + z[\partial^2 \times \ln f(z; \boldsymbol{\tau}_0)/\partial z^2]| < \infty$, the absolute value of (A.46) is bounded above by

$$(\text{A.47}) \qquad (constant) \sup_{i \in \{1, \ldots, p\}} \left(\frac{|\hat{\theta}_{i,\mathrm{ML}} - \theta_{0i}|}{n-p} \sum_{t=p+1}^{n} \left|\frac{\partial Z_t(\boldsymbol{\theta}_n^*, s_0)}{\partial \theta_i}\right|\right).$$

Recall, from the Proof of Lemma A.5, that $\kappa_2 = \alpha_0(1 + \alpha_0/3)/(1 + \alpha_0/2)I\{\alpha_0 \leq 1\} + I\{\alpha_0 > 1\}$. For $i \in \{1, \ldots, p\}$ and $\epsilon > 0$,

$$\mathrm{P}\left(\left[\frac{1}{(n-p)^{1/2+1/\alpha_0}} \sum_{t=p+1}^{n} \left|\frac{\partial Z_t(\boldsymbol{\theta}_n^*, s_0)}{\partial \theta_i}\right| I\{\|\hat{\boldsymbol{\theta}}_{\mathrm{ML}} - \boldsymbol{\theta}_0\| < \delta\}\right]^{\kappa_2} > \epsilon^{\kappa_2}\right)$$

$$\leq \epsilon^{-\kappa_2}(n-p)^{1-\kappa_2(1/2+1/\alpha_0)} \mathrm{E}\left\{\left|\frac{\partial Z_t(\boldsymbol{\theta}_n^*, s_0)}{\partial \theta_i}\right|^{\kappa_2} I\{\|\hat{\boldsymbol{\theta}}_{\mathrm{ML}} - \boldsymbol{\theta}_0\| < \delta\}\right\},$$

which can be shown to be $o(1)$ for sufficiently small $\delta > 0$ since $\kappa_2(1/2 + 1/\alpha_0) > 1$ and $\mathrm{E}|Z_1|^{\kappa_2} < \infty$. Therefore, since $n^{1/\alpha_0}(\hat{\boldsymbol{\theta}}_{\mathrm{ML}} - \boldsymbol{\theta}_0) = O_p(1)$, it follows that (A.47), and hence (A.45) and (A.46), are $o_p(n^{-1/2})$, and so (A.43) is $O_p(n^{-1/2})$. Since $m_n/n \to 0$, (A.42) must be $o_p(1)$, and so the proof is complete. $\quad\square$

**Acknowledgments.** We wish to thank two anonymous reviewers for their helpful comments.

B. ANDREWS
DEPARTMENT OF STATISTICS
NORTHWESTERN UNIVERSITY
2006 SHERIDAN ROAD
EVANSTON, ILLINOIS 60208
USA
E-MAIL: bandrews@northwestern.edu

M. CALDER
PHZ CAPITAL PARTNERS
321 COMMONWEALTH ROAD
WAYLAND, MASSACHUSETTS 01778
USA
E-MAIL: calder@phz.com

R. A. Davis
Department of Statistics
Columbia University
1255 Amsterdam Avenue
New York, New York 10027
USA
E-mail: rdavis@stat.columbia.edu